\documentclass[a4paper]{article}

\usepackage[applemac]{inputenc}
\usepackage[T1]{fontenc}
\usepackage{lmodern}
\usepackage[english]{babel}
\usepackage{amssymb}
\usepackage{amsmath}
\usepackage{theorem}
\usepackage{graphicx}
\usepackage{bbm}

\newcommand{\R}{\mathbb{R}}
\newcommand{\N}{\mathbb{N}}
\newcommand{\PP}{\mathbf{P}}
\newcommand{\E}{\mathbf{E}}
\newcommand{\1}{\mathbbm{1}}
\newcommand{\demo}{\paragraph{Proof}}
\newcommand{\findemo}{\hfill$\Box$}
\newcommand{\Lc}{\mathbf{L}_{c}}
\newcommand{\tstar}{t^{*}}
\newcommand{\ts}{t^{*}}
\def\union{\mathop{\cup}}

\newtheorem{theorem}{Theorem}[section]
\newtheorem{lemme}[theorem]{Lemma}
\newtheorem{proposition}[theorem]{Proposition}
\newtheorem{definition}[theorem]{Definition}
\newtheorem{hyp}[theorem]{Assumption}
\newtheorem{remarque}[theorem]{Remark}

\title{Numerical method for expectations of piecewise-determistic Markov processes
\thanks{This work was supported by ARPEGE program of the French National Agency of Research (ANR),
project ``FAUTOCOES'', number ANR-09-SEGI-004.}
}

\author{ \mbox{ }
A. Brandejsky
\footnote{Postal address: IMB, Universit\'e Bordeaux 1, 351 cours de la lib\'eration, 33405 Talence cedex, France.}
\\ 
\small INRIA Bordeaux Sud-Ouest, team CQFD, F-33400 Talence, France.\\
\small Univ. Bordeaux, IMB, UMR 5251, F-33400 Talence, France. \\
\small CNRS, IMB, UMR 5251, F-33400 Talence, France. 
\and
B. de Saporta $^{\dag}$
\\
\small Univ. Bordeaux, Gretha, UMR 5113, F-33400 Talence, France. \\
\small CNRS, Gretha, UMR 5113, F-33400 Talence, France. \\
\small CNRS, IMB, UMR 5251, F-33400 Talence, France. \\
\small INRIA Bordeaux Sud-Ouest, team CQFD, F-33400 Talence, France. 
\\
\and
F. Dufour $^{\dag}$
\\
\small Univ. Bordeaux, IMB, UMR 5251, F-33400 Talence, France. \\
\small CNRS, IMB, UMR 5251, F-33400 Talence, France. \\
\small INRIA Bordeaux Sud-Ouest, team CQFD, F-33400 Talence, France. 
}

\date\today

\begin{document}
\maketitle

\vspace{2cm}
\paragraph{Abstract} We present a numerical method to compute expectations of functionals of a piecewise-deterministic Markov process. We discuss time dependent functionals as well as deterministic time horizon problems. Our approach is based on the quantization of an underlying discrete-time Markov chain. We obtain bounds for the rate of convergence of the algorithm. The approximation we propose is easily computable and is flexible with respect to some of the parameters defining the problem. Two examples illustrate the paper.


\paragraph{Maths Subject Classification 2010} Primary: 60J25, 65C20. Secondary: 60K10.

\newpage

\tableofcontents

\newpage

\section{Introduction}

The aim of this paper is to propose a practical numerical method to approximate some expectations related to a piecewise-deterministic Markov process thanks to the quantization of a discrete-time Markov chain naturally embedded within the continuous-time process.\\

Piecewise-deterministic Markov processes (PDMP's) have been introduced by M.H.A. Davis in \cite{davis93} as a general class of stochastic models. PDMP's are a family of Markov processes involving deterministic motion punctuated by random jumps. The motion depends on three local characteristics namely the flow $\Phi$,  the jump rate $\lambda$ and the transition measure $Q$,which specifies the post-jump location. Starting from the point $x$, the motion of the process follows the flow $\Phi(x,t)$ until the first jump time $T_{1}$, which occurs either spontaneously in a Poisson-like fashion with rate $\lambda(\Phi(x,t))$ or when the flow $\Phi(x,t)$ hits the boundary of the state space. In either case, the location of the process at the jump time $T_{1}$, is selected by the transition measure $Q(\Phi(x,T_{1}),\cdot)$ and the motion restarts from this new point $X_{T_{1}}$ denoted $Z_{1}$. We define similarly the time $S_{2}$ until the next jump, $T_{2}=T_{1}+S_{2}$ with the next post-jump location defined by $Z_{2}=X_{T_{2}}$ and so on. Thus, associated to the PDMP we have the discrete-time Markov chain $(Z_{n},S_{n})_{n\in \N}$, given by the post-jump locations and the inter-jump times. A suitable choice of the state space and the local characteristics $\Phi$, $\lambda$ and $Q$ provides stochastic models covering a great number of problems of operations research as described in \cite{davis93} section 33.\\

We are interested in the approximation of expectations of the form
$$\E_{x}\left[\int_{0}^{T_{N}}l(X_{t})dt+\sum_{j=1}^{N}c(X_{T^{-}_{j}})\1_{\{X_{T^{-}_{j}}\in \partial E\}}\right]$$
where $\left(X_{t}\right)_{t\geq 0}$ is a PDMP and $l$ and $c$ are some non negative, real-valued, bounded functions and $\partial E$ is the boundary of the domain. Such expectations are discussed by M.H.A. Davis in \cite{davis93}, chapter 3. They often appear as ``cost'' or ``reward'' functions in optimization problems. The first term is referred to as the running cost while the second may be called the boundary jump cost. Besides, they are quite general since M.H.A. Davis shows how a ``wide variety of apparently different functionals'' can be obtained from the above specific form. For example, this wide variety includes quantities such as a mean exit time and even, for any fixed $t\geq 0$, the distribution of $X_{t}$ (i.e. $\E_{x}[\1_{F}(X_{t})]$ where $F$ is a measurable set).\\

There are surprisingly few works in the literature devoted to the actual computation of such expectations, using other means than direct Monte Carlo simulations. M.H.A Davis showed that these expectations satisfy integro-differential equations. However, the set of partial differential equations that is obtained is unusual. Roughly speaking, these differential equations are basically transport equations with a non-constant velocity and they are coupled by the boundary conditions and by some integral terms involving kernels that are derived from the properties of the underlying stochastic process. The main difficulty comes from the fact that the domains on which the equations have to be solved vary from one equation to another making their numerical resolution highly problem specific. Another similar approach has been recently investigated in \cite{cocozza06, eymard08}. It is based on a discretization of the Chapman Kolmogorov equations satisfied by the distribution of the process $\left(X_{t}\right)_{t\geq 0}$. The authors propose an approximation of such expectations based on finite volume methods. Unfortunately, their method is only valid if there are no jumps at the boundary. Our approach is completely different and does not rely on differential equations, but on the fact that such expectations can be computed by iterating an integral operator $G$. This operator only involves the embedded Markov chain  $(Z_{n},S_{n})_{n\in \N}$ and conditional expectations. It is therefore natural to propose a computational method based on the quantization of this Markov chain, following the same idea as \cite{saporta10}.\\

There exists an extensive literature on quantization methods for random variables and processes. The interested reader may for instance consult \cite{gray98}, \cite{pages04} and the references within. Quantization methods have been developed recently in numerical probability or optimal stochastic control with applications in finance (see e.g. \cite{bally03}, \cite{bally05} and \cite{pages04}). The quantization of a random variable $X$ consists in finding a finite grid such that the projection $\widehat{X}$ of $X$ on this grid minimizes some $L^{p}$ norm of the difference $X-\widehat{X}$. Roughly speaking, such a grid will have more points in the areas of high density of $X$. As explained for instance in \cite{pages04}, section 3, under some Lipschitz-continuity conditions, bounds for the rate of convergence of functionals of the quantized process towards the original process are available.\\

In the present work, we develop a numerical method to compute expectations of functionals of the above form where the cost functions $l$ and $c$  satisfy some Lipschitz-continuity conditions. We first recall the results presented by M.H.A. Davis according to whom, the above expectation may be computed by iterating an operator denoted $G$. Consequently, it appears natural to follow the idea developed in \cite{saporta10} namely to express the operator $G$ in terms of the underlying discrete-time Markov chain $(Z_{n},S_{n})_{n\in \N}$ and to replace it by its quantized approximation. Moreover, in order to prove the convergence of our algorithm, we replace the indicator function $\1_{\{X_{T^{-}_{j}}\in \partial E\}}$ contained within the functional by some Lipschitz continuous approximation. Bounds for the rate of convergence are then obtained. However, and this is the main contribution of this paper, we then tackle two important aspects that had not been investigated in \cite{saporta10}.\\

The first aspect consists in allowing $c$ and $l$ to be time depending functions, although still Lipschitz continuous, so that we may compute expectations of the form 
$$\E_{x}\left[\int_{0}^{T_{N}}l(X_{t},t)dt+\sum_{j=1}^{N}c(X_{T^{-}_{j}},T_{j})\1_{\{X_{T^{-}_{j}}\in \partial E\}}\right].$$
This important generalization has huge applicative consequences. For instance, it allows discounted ``cost'' or ``reward'' functions such as $l(x,t)=e^{-\delta t}l(x)$ and $c(x,t)=e^{-\delta t}c(x)$ where $\delta$ is some interest rate. To compute the above expectation, our strategy consists in considering, as it is suggested by M.H.A. Davis in \cite{davis93}, the time augmented process~$\widetilde{X}_{t}=(X_{t},t)$. Therefore, a natural way to deal with the time depending problem is to apply our previous approximation scheme to the time augmented process $(\widetilde{X}_{t})_{t\geq 0}$. However, it is far from obvious, that the assumptions required by our numerical method still hold for this new PDMP $(\widetilde{X}_{t})_{t\geq 0}$.\\

The second important generalization is to consider the deterministic time horizon problem. Indeed, it seems crucial, regarding the applications, to be able to approximate 
\begin{align*}
\E_{x}\Big[\int_{0}^{t_{f}}l(X_{t},t)dt+\sum_{T_{j}\leq t_{j}}c(X_{T^{-}_{j}},T_{j})\1_{\{X_{T^{-}_{j}}\in \partial E\}}\Big]
\end{align*}
 for some fixed $t_{f}>0$ regardless of how many jumps occur before this deterministic time. To compute this quantity, we start by choosing a time~$N$ such that $\PP\left(T_{N}<t_{f}\right)$ be small so that the previous expectation boils down to $\E_{x}\left[\int_{0}^{T_{N}}l(X_{t},t)\1_{\{t\leq t_{f}\}}dt+\sum_{j=1}^{N}c(X_{T^{-}_{j}},T_{j})\1_{\{X_{T^{-}_{j}}\in \partial E\}}\1_{\{T_{j}\leq t_{f}\}}\right]$. At first sight, this functional seems to be of the previous form. Yet, one must recall that Lipschitz continuity conditions have been made concerning the cost functions so that the indicator functions $\1_{\{\cdot\leq t_{f}\}}$ prevent a direct application of the earlier results. We deal with the two indicator functions in two different ways. On the one hand, we prove that it is possible to relax the regularity condition on the running cost function so that our algorithm still converges in spite of the first indicator function. On the other hand, since the same reasoning cannot be applied to the indicator function within the boundary jump cost term, we bound it between two Lipschitz continuous functions. This provides bounds for the expectation of the deterministic time horizon functional.\\

An important advantage of our method is that it is flexible. Indeed, as pointed out in \cite{bally03}, a quantization based method is ``obstacle free'' which means, in our case, that it produces, once and for all, a discretization of the process independently of the functions $l$ and $c$ since the quantization grids merely depend on the dynamics of the process. They are only computed once, stored off-line and may therefore serve many purposes. Once they have been obtained, we are able to approximate very easily and quickly any of the expectations described earlier. This flexibility is definitely an important advantage of our scheme over standard methods such as Monte-Carlo simulations since, with such methods, we would have to run the whole algorithm for each expectation we want to compute. This point is illustrated in Section~\ref{section-result} where we easily solve an optimization problem that would be very laboriously handled by Monte-Carlo simulations.\\

The paper is organized as follows. We first recall, in Section~\ref{sec-def}, the definition of a PDMP and state our assumptions. In Section~\ref{section-exit-time}, we introduce the recursive method to compute the expectation. Section~\ref{section-approximation} presents the approximation scheme and a bound for the rate of convergence. The main contribution of the paper lies in Section~\ref{section-transformation} which contains the generalizations to the time dependent parameters and the deterministic time horizon problems. Eventually, the paper is illustrated by two numerical examples in Section~\ref{section-result} and concluded in Section~\ref{section-conclusion} while technical results are postponed to the Appendix.

\section{Definitions and assumptions}\label{sec-def}

For all metric space $E$, we denote $\mathcal{B}(E)$ its Borel $\sigma$-field and $B(E)$ the set of real-valued, bounded and measurable functions defined on~$E$. For $a,b \in \R$, denote $a\wedge b=\min(a,b)$, $a\vee b=\max(a,b)$ and $a^{+}=a \vee 0$.

\subsection*{Definition of a PDMP}
In this first section, let us define a piecewise-deterministic Markov process and introduce some general assumptions. Let $M$ be a finite set called the set of the modes that will represent the different regimes of evolution of the PDMP. For each $m\in M$, the process evolves in $E_m$, an open subset of $\R^{d}$. Let $$E=\left\{(m,\zeta),m\in M,\zeta\in E_m \right\}.$$
This is the state space of the process $(X_t)_{t\in \mathbb{R}^+}=(m_t,\zeta_t)_{t\in \mathbb{R}^+}$. Let $\partial E$ be its boundary and $\overline{E}$ its closure and for any subset $Y$ of $E$, $Y^{c}$ denotes its complement.
\paragraph{}
A PDMP is defined by its local characteristics $(\Phi_{m},\lambda_{m},Q_{m})_{m\in M}$.
\begin{itemize}
\item{For each $m\in M$, $\Phi_m : \mathbb{R}^{d}\times \mathbb{R}\rightarrow \mathbb{R}^{d}$ is a continuous function called the flow in mode $m$. For all $t\in \mathbb{R}$, $\Phi_m(\cdot,t)$ is an homeomorphism and $t\rightarrow \Phi_m(\cdot,t)$ is a semi-group i.e. for all $\zeta\in \mathbb{R}^{d}$, $\Phi_m(\zeta,t+s)=\Phi_m(\Phi_m(\zeta,s),t)$.
    For all $x=(m,\zeta)\in E$, define now the deterministic exit time from $E$ :
    $$t^*(x)=
    \inf\{t>0 \text{ such that } \Phi_m(\zeta,t)\in\partial E_m\}.$$
    We use here and throughout the whole paper the convention $\inf \emptyset = + \infty$.
    }
\item{For all $m\in M$, the jump rate $\lambda_{m} : \overline{E}_{m}\rightarrow \mathbb{R}^+$ is measurable and satisfies~:
$$\forall (m,\zeta)\in E\text{,  }\exists \epsilon >0\text{ such that }\int_0^\epsilon \lambda_{m}(\Phi_m(\zeta,t))dt< +\infty.$$
}
\item{For all $m\in M$, $Q_{m}$ is a Markov kernel on $(\mathcal{B}(\overline{E}),\overline{E}_{m})$ which satisfies :
$$\forall \zeta \in \overline{E}_{m}\text{,  }Q_{m}(\zeta,\{(m,\zeta)\}^{c})=1.$$
}
\end{itemize}
From these characteristics, it can be shown (see \cite{davis93}) that there exists a filtered probability space $(\Omega,\mathcal{F},{\mathcal{F}_t},(\mathbf{P}_x)_{x\in E})$ on which a process $(X_t)_{t\in \mathbb{R}^+}$ is defined. Its motion, starting from a point $x\in E$, may be constructed as follows.
Let $T_1$ be a nonnegative random variable with survival function :
$$\PP_x(T_1>t)=
\left\{\begin{array}{ll}
e^{-\Lambda(x,t)} & \text{if } 0\leq t<t^*(x), \\
0 & \text{if }t\geq t^*(x),
\end{array}\right.$$
where for $x=(m,\zeta)\in E$ and $t\in [0,t^*(x)]$,
$$\Lambda(x,t)=\int_0^t\lambda_{m}(\Phi_m(\zeta,s))ds.$$
One then chooses an $E$-valued random variable $Z_1$ according to the distribution $Q_{m}(\Phi_m(\zeta,T_1),\cdot)$. The trajectory of $X_t$ for $t\leq T_1$ is :
$$X_t=\left\{
\begin{array}{ll}
(m,\Phi_m(\zeta,t))&\text{ if }t<T_1,\\
Z_1&\text{ if }t=T_1.
\end{array}\right.
$$
Starting from the point $X_{T_1}=Z_1$, one then selects in a similar way $S_{2}=T_2-T_1$ the time between $T_{1}$ and the next jump, $Z_2$ the next post-jump location and so on. M.H.A. Davis shows, in \cite{davis93}, that the process so defined is a strong Markov process $(X_t)_{t\geq 0}$ with jump times $(T_n)_{n\in {\mathbb{N}}}$ (with $T_0=0$). The process $(\Theta_n)_{n\in\N}=(Z_n,S_n)_{n\in\N}$ where $Z_n=X_{T_n}$ is the post-jump location and $S_n=T_n-T_{n-1}$ (with $S_0=0$) is the $n$-th inter-jump time is clearly a discrete-time Markov chain. \\
The following assumption about the jump-times is standard (see for example \cite{davis93}, section 24) :
\begin{hyp}\label{Tk_goes_to_infty}
For all $(x,t)\in E\times \mathbb{R}^+$, $\E_x\left[\sum_k \mathbbm{1}_{\{T_k<t\}}\right]<+\infty$.
\end{hyp}
It implies in particular that $T_k$ goes to infinity a.s. when $k$ goes to infinity.

\subsection*{Notation and assumptions}

For notational convenience, any function $h$ defined on $E$ will be identified with its component functions $h_m$ defined on $E_m$. Thus, one may write $$h(x)=h_m(\zeta) \text{ when } x=(m,\zeta)\in E.$$
We also define a generalized flow $\Phi : E\times\mathbb{R}^+\rightarrow E$ such that $$\Phi(x,t)=(m,\Phi_m(\zeta,t)) \text{ when } x=(m,\zeta)\in E.$$

\paragraph{}Define on $E$ the following distance, for $x=(m,\zeta)$ and $x'=(m',\zeta') \in E$,
\begin{equation}\label{def-dist}
|x-x'|=\left\{
\begin{array}{ll}
+\infty & \text{ if } m\neq m',\\
|\zeta-\zeta'|& \text{ otherwise. }
\end{array}
\right.
\end{equation}
For any function $w$ in $B(\overline{E})$, introduce the following notation 
\begin{equation*}
Qw(x)=\int_E w(y)Q(x,dy) \text{, }\qquad C_w=\sup_{x\in \overline{E}}|w(x)|,
\end{equation*}
and for any Lipschitz continuous function $w$ in $B(E)$, denote $[w]^{E}$, or if there is no ambiguity $[w]$, its Lipschitz constant:
$$[w]^{E}=\sup_{x\neq y \in E}\frac{|w(x)-w(y)|}{|x-y|},$$
with the convention $\frac{1}{\infty}=0$.

\begin{remarque}\label{RqLip_f_fm}
For $w\in B(\overline{E})$ and from the definition of the distance on $E$, one has $[w]=\text{max}_{m\in M}[w_m]$.
\end{remarque}

\begin{definition}\label{def-Lc}
Denote $\Lc(E)$ the set of functions $w\in B(E)$ that are Lipschitz continuous along the flow i.e. the real-valued, bounded, measurable functions defined on $E$ and satisfying the following conditions:
\begin{itemize}
\item{For all $x\in E$, $w(\Phi(x,\cdot)) \text{ : } [0,t^{*}(x)[\rightarrow \R$ is continuous, $\lim_{t\rightarrow t^{*}(x)}w(\Phi(x,t))$ exists and is denoted $w\big(\Phi(x,t^{*}(x))\big)$,}
\item{there exists $[w]^{E}_1\in \mathbb{R}^+$ such that for all $x,y\in E$ and $t\in [0,t^*(x)\wedge t^*(y)]$, one has:
$$|w(\Phi(x,t))-w(\Phi(y,t))|\leq [w]^{E}_1|x-y|,$$
}
\item{there exists $[w]^{E}_2\in \mathbb{R}^+$ such that for all $x\in E$ and $t,u\in [0,t^*(x)]$, one has:
$$|w(\Phi(x,t))-w(\Phi(x,u))|\leq [w]^{E}_2|t-u|,$$}
\item{there exists $[w]^{E}_*\in \mathbb{R}^+$ such that for all $x,y\in E$, one has:
$$|w(\Phi(x,t^*(x)))-w(\Phi(y,t^*(y)))|\leq [w]^{E}_*|x-y|.$$}
\end{itemize}
Denote also $\Lc(\partial E)$ the set of real-valued, bounded, measurable functions defined on $\partial E$ satisfying the following condition:
\begin{itemize}
\item{there exists $[w]^{\partial E}_*\in \mathbb{R}^+$ such that for all $x,y\in E$, one has:
$$|w(\Phi(x,t^*(x)))-w(\Phi(y,t^*(y)))|\leq [w]^{\partial E}_*|x-y|.$$}
\end{itemize}
\end{definition}

\begin{remarque} When there is no ambiguity, we will denote $[w]_i$ instead of $[w]^{E}_i$ for $i\in\{1,2,*\}$ and $[w]_*$ instead of $[w]^{\partial E}_*$.
\end{remarque}

\begin{remarque}
In the above definition, we used the generalized flow for notational convenience. For instance, the definition of $[w]_1$ is equivalent to the following: for all $m\in M$, there exists $[w_m]_1\in \mathbb{R}^+$ such that for all $\zeta,\zeta'\in E_m$ and $t\in [0,t^*(m,\zeta)\wedge t^*(m,\zeta')]$, one has:
$$|w_m(\Phi_m(\zeta,t))-w_m(\Phi_m(\zeta',t))|\leq [w_m]_1|\zeta-\zeta'|.$$
Let $[w]_1=\text{max }_{m\in M}[w_m]_1$.
\end{remarque}

\begin{definition}\label{def-Lcu}
For all $u\geq 0$, denote $\Lc^{u}(E)$ the set of functions $w\in B(E)$ Lipschitz continuous along the flow until time $u$ i.e. the real-valued, bounded, measurable functions defined on $E$ and satisfying the following conditions:
\begin{itemize}
\item{For all $x\in E$, $w(\Phi(x,\cdot)) \text{ : } [0,t^{*}(x)\wedge u[\rightarrow \R$ is continuous and if $\ts(x)\leq u$, then $\lim_{t\rightarrow t^{*}(x)}w(\Phi(x,t))$ exists and is denoted $w\big(\Phi(x,t^{*}(x))\big)$,}
\item{there exists $[w]^{E,u}_1\in \mathbb{R}^+$ such that for all $x,y\in E$ and $t\in [0,t^*(x)\wedge t^*(y)\wedge u]$, one has:
$$|w(\Phi(x,t))-w(\Phi(y,t))|\leq [w]^{E,u}_1|x-y|,$$
}
\item{there exists $[w]^{E,u}_2\in \mathbb{R}^+$ such that for all $x\in E$ and $t,t'\in [0,t^*(x)\wedge u]$, one has:
$$|w(\Phi(x,t))-w(\Phi(x,t'))|\leq [w]^{E,u}_2|t-t'|,$$}
\item{there exists $[w]^{E,u}_*\in \mathbb{R}^+$ such that for all $x,y\in E$, if $\ts(x)\leq u$ and $\ts(y)\leq u$, one has:
$$|w(\Phi(x,t^*(x)))-w(\Phi(y,t^*(y)))|\leq [w]^{E,u}_*|x-y|.$$}
\end{itemize}
\end{definition}

\begin{remarque}For all $u\leq u'$, one has $\Lc^{u'}(E)\subset \Lc^{u}(E)$ with $[w]^{E,u}_i\leq [w]^{E,u'}_i$ where $i\in \{1,2,*\}$.
\end{remarque}

\begin{remarque}Note that Definitions \ref{def-Lc} and \ref{def-Lcu} correspond respectively to the Lipschitz and local Lipschitz continuity along the flow that is, along the trajectories of the process. They can be replaced by (local) Lipschitz assumptions on the flow $\Phi$, $\ts$ and $w$ in the classical sense.
\end{remarque}

We will require the following assumptions.

\begin{hyp}\label{hyp-lip-lambda}
The jump rate $\lambda$ is bounded and there exists $[\lambda]_1\in \mathbb{R}^+$ such that for all $x,y\in E$ and $t\in [0,t^*(x)\wedge t^*(y)]$, one has:
$$|\lambda(\Phi(x,t))-\lambda(\Phi(y,t))|\leq [\lambda]_1|x-y|.$$
\end{hyp}

\begin{hyp}\label{hyp-lip-t*}
The deterministic exit time from $E$, denoted $t^*$, is assumed to be bounded and Lipschitz continuous on $E$.
\end{hyp}

\begin{remarque}Since the deterministic exit time $\ts$ is bounded by $C_{\ts}$, one may notice that $\Lc^{u}(E)$ for $u\geq C_{\ts}$ is no other than $\Lc(E)$.
\end{remarque}

\begin{remarque} In most practical applications, the physical properties of the system ensure that either $\ts$ is bounded, or the problem has a natural finite deterministic time horizon $t_{f}$. In the latter case, there is no loss of generality in considering that $\ts$ is bounded by this deterministic time horizon. This leads to replacing~$C_{\ts}$ by $t_{f}$. An example of such a situation is presented in an industrial example in Section \ref{exemple-corrosion}.
\end{remarque}
 
\begin{hyp}\label{hyp-lip-Q}
The Markov kernel $Q$ is Lipschitz in the following sense: there exists $[Q]\in \mathbb{R}^+$ such that for all $u\geq 0$ and for all function $w\in \Lc^{u}(E)$, one has:
\begin{enumerate}
\item{for all $x,y\in E$ and $t\in [0,t^*(x)\wedge t^*(y)\wedge u[$,
$$|Qw(\Phi(x,t))-Qw(\Phi(y,t))|\leq [Q][w]^{E,u}_1|x-y|.$$}
\item{for all $x,y\in E$ such that $\ts(x)\vee\ts(y)\leq u$,
$$|Qw(\Phi(x,t^*(x)))-Qw(\Phi(y,t^*(y)))|\leq [Q]\big([w]^{E,u}_*+[w]^{E,u}_1\big)|x-y|.$$}
\end{enumerate}
\end{hyp}

\begin{remarque} Notice that assumption~\ref{hyp-lip-Q} is slightly more restrictive that its counterpart in \cite{saporta10} (assumption 2.5) because of the introduction of the state space $\Lc^{u}(E)$. This is to ensure that the time augmented process still satisfies a similar assumption, see Section~\ref{section-time-augemented}.

\end{remarque}

\section{Expectation}\label{section-exit-time}\label{DefRec}

From now on, we will assume that $Z_{0}=x$ a.s. for some $x\in E$. For all fixed~$N\in\N^{*}$, we intend to numerically approximate the quantity
\begin{equation}\label{J}
J_{N}(l,c)(x)=\E_{x}\left[\int_{0}^{T_{N}}l(X_{t})dt+\sum_{j=1}^{N}c(X_{T^{-}_{j}})\1_{\{X_{T^{-}_{j}}\in \partial E\}}\right],
\end{equation}
where $l\in B(E)$, $c\in B(\partial E)$ and $X_{t^{-}}$ is the left limit of $X_{t}$. Thus, $X_{T^{-}_{j}}$ is the $j$-th pre-jump location. Since the boundary jumps occur exactly at the deterministic exit times from $E$, one has,
\begin{align*}
J_{N}(l,c)(x)&=\E_{x}\left[\int_{0}^{T_{N}}l(X_{t})dt+\sum_{j=1}^{N}c\big(\Phi(Z_{j-1},\tstar(Z_{j-1}))\big)\1_{\{S_{j}=\tstar(Z_{j-1})\}}\right]
\end{align*}
In many applications, $J_{N}(l,c)(x)$ appears as a ``cost'' or a ``reward'' function. The first term, that depends on $l$, is called the running cost and the second one, that depends on $c$, is the boundary jump cost.

The rest of this section is dedicated to finding a formulation of the above expectation that will allow us to derive a numerical computation method. The Lipschitz continuity property will be a crucial point when it will come to proving the convergence of our approximation scheme. For this reason, the first step of our approximation is to replace the indicator function in $J_{N}(l,c)(x)$ by a Lipschitz continuous function. Then, we will present a recursive method yielding the required expectation. This recursive formulation will be the basis of our numerical method.

\subsection{Lipschitz continuity}

We introduce a regularity assumption on $l$ and $c$.
\begin{hyp}\label{hyp-lip-lc}
We assume that $l\in\Lc(E)$ and $c\in\Lc(\partial E)$.
\end{hyp}

Moreover, we replace the indicator function in $J_{N}(l,c)(x)$ by a Lipschitz continuous function denoted $\delta^{A}$ with $A>0$. Let then
$$J_{N}^{A}(l,c)(x)=\E_{x}\left[\int_{0}^{T_{N}}l(X_{t})dt+\sum_{j=1}^{N}c\big(\Phi(Z_{j-1},\tstar(Z_{j-1}))\big)\delta^{A}(Z_{j-1},S_{j})\right],$$
where $\delta^{A}$ is a triangular approximation of the indicator function. It is defined on $E\times\R$ by
$$\delta^{A}(x,t)=\left\{\begin{array}{ll}
A\big(t-(\tstar(x)-\frac{1}{A})\big) & \text{ for } t\in [\tstar(x)-\frac{1}{A};\tstar(x)],\\
-A\big(t-(\tstar(x)+\frac{1}{A})\big) & \text{ for } t\in [\tstar(x);\tstar(x)+\frac{1}{A}],\\
0 & \text{ otherwise}.
\end{array}\right.$$

For all $x\in E$, the function $\delta^{A}(x,t)$ goes to $\1_{\{t=\tstar(x)\}}$ when $A$ goes to infinity. The following proposition proves the convergence of $J_{N}^{A}(l,c)(x)$ towards $J_{N}(l,c)(x)$ with an error bound.

\begin{proposition}\label{prop-JAconv} For all $x\in E$, $A>0$, $N\in\N^{*}$, $l\in\Lc(E)$ and $c\in\Lc(\partial E)$, one has
$$\left|J_{N}^{A}(l,c)(x)-J_{N}(l,c)(x)\right|\leq\frac{NC_{c}C_{\lambda}}{A}.$$
\end{proposition}

\demo For all $x\in E$, one has
\begin{eqnarray*}\lefteqn{\left|J_{N}^{A}(l,c)(x)-J_{N}(l,c)(x)\right|}\\
&=&\left|\E_{x}\left[\sum_{j=1}^{N}c\big(\Phi(Z_{j-1},\tstar(Z_{j-1}))\big)\big(\delta^{A}(Z_{j-1},S_{j})-\1_{\{S_{j}=\tstar(Z_{j-1})\}}\big)\right]\right|\\
&\leq& C_{c}\sum_{j=1}^{N}\E_{x}\left[\big|\delta^{A}(Z_{j-1},S_{j})-\1_{\{S_{j}=\tstar(Z_{j-1})\}}\big|\right]\\
&\leq& C_{c}\sum_{j=1}^{N}\E_{x}\left[\E\left[\big|\delta^{A}(Z_{j-1},S_{j})-\1_{\{S_{j}=\tstar(Z_{j-1})\}}\big|\Big|Z_{j-1}\right]\right].
\end{eqnarray*}
We recall that the conditional law of $S_{j}$ with respect to $Z_{j-1}$ has density $s\rightarrow \lambda\big(\Phi(Z_{j-1},s)\big)e^{-\Lambda(Z_{j-1},s)}$ on $[0;\ts(Z_{j-1})[$ and puts the weight $e^{-\Lambda(Z_{j-1},\ts(Z_{j-1}))}$ on the point $\ts(Z_{j-1})$. We also recall that $\lambda$ is bounded thanks to assumption~\ref{hyp-lip-lambda}. Eventually, one has
\begin{eqnarray*}
\lefteqn{\left|J_{N}^{A}(l,c)(x)-J_{N}(l,c)(x)\right|}\\
&\leq& C_{c}\sum_{j=1}^{N}\E_{x}\left[\int_{\tstar(Z_{j-1})-\frac{1}{A}}^{\tstar(Z_{j-1})}\delta^{A}(Z_{j-1},s)\lambda\big(\Phi(Z_{j-1},s)\big)e^{-\Lambda(Z_{j-1},s)}ds\right]\\
&\leq &\frac{NC_{c}C_{\lambda}}{A}.
\end{eqnarray*}
Hence the result.
\findemo\\

Consequently to this proposition, we consider, from now on, the approximation of $J_{N}^{A}(l,c)(x)$ for some fixed $A$, large enough to ensure that the previous error is as small as required. The suitable choice of $A$ will be discussed in section~\ref{choixA}.

\subsection{Recursive formulation}

M.H.A. Davis shows in \cite{davis93}, section 32, that the expectation $J_{N}^{A}(l,c)(x)$ we are interested in is obtained by merely iterating an operator that we will denote $G$. The rest of this section is dedicated to presenting this method from which we will derive our approximation scheme in Section~\ref{section-approximation}.

\begin{definition}\label{def-L-C-F-G}Introduce the functions $L$, $C$ and $F$ defined for all $x\in E$ and $t\in[0;t^{*}(x)]$ by
\begin{align*}
L(x,t)&=\int_{0}^{t}l\big(\Phi(x,s)\big)ds,\\
C(x,t)&=c\big(\Phi(x,\tstar(x))\big)\delta^{A}(x,t),\\
F(x,t)&=L(x,t)+C(x,t),
\end{align*}
along with the operator $G$ : $B(E)\rightarrow B(E)$
$$Gw(x)=\E_{x}\left[F(x,S_{1})+w(Z_{1})\right].$$
\end{definition}

\begin{definition}
Introduce the sequence of functions $(v_{k})_{0\leq k\leq N}$ in $B(E)$ defined as follows:
$$\left\{\begin{array}{ll}
v_{N}(x)&=0,\\
v_{k}(x)&=Gv_{k+1}(x),
\end{array}\right.$$
\end{definition}
M.H.A. Davis then shows in \cite{davis93}, equation 32.33, that for all $k\in \{0,...,N\}$,
$$v_{N-k}(x)=\E_{x}\left[\int_{0}^{T_{k}}l(X_{t})dt+\sum_{j=1}^{k}c\big(\Phi(Z_{j-1},\tstar(Z_{j-1}))\big)\delta^{A}(Z_{j-1},S_{j})\right].$$
Thus, the quantity $J_{N}^{A}(l,c)(x)$ we intend to approximate is none other than~$v_{0}(x)$.\\

Notice that, thanks to the Markov property of the chain $(Z_{n},S_{n})_{n\in\N}$, one has for all $k\in \{0,...,N-1\}$,
\begin{equation}\label{markov-op-G}
Gw(x)=\E\left[F(Z_{k},S_{k+1})+w(Z_{k+1})\big|Z_{k}=x\right].
\end{equation}
Hence, for all $k\in \{0,...,N\}$, let $V_{k}=v_{k}(Z_{k})$ so that one has
$$\left\{\begin{array}{ll}
V_{N}&=0,\\
V_{k}&=\E\left[F(Z_{k},S_{k+1})+V_{k+1}\big|Z_{k}\right].
\end{array}\right.$$
This backward recursion provides the required quantity
$$V_{0}=J_{N}^{A}(l,c)(x).$$

Consequently, we need to approximate the sequence of random variables $(V_{k})_{0\leq k\leq N}$. This sequence satisfies a recursion that only depends on the chain $(Z_{k},S_{k})_{0\leq k\leq N}$. Therefore, it appears natural to propose an approximation scheme based on a discretization of this chain $(Z_{k},S_{k})_{0\leq k\leq N}$, called quantization, similarly to the ideas developed in \cite{saporta10} and \cite{brandejsky10}.

\section{Approximation scheme}\label{section-approximation}

Let us now turn to the approximation scheme itself. We explained in the previous section how the expectation we are interested in stems from the iteration of the operator $G$ that only depends on the discrete-time Markov chain $(Z_{k},S_{k})_{0\leq k\leq N}$. The first step of our numerical method is therefore to discretize this chain in order to approximate the operator $G$.

\subsection{Quantization of the chain $(Z_k,S_k)_{k\leq N}$}

Our approximation method is based on the quantization of the underlying discrete time Markov chain $(\Theta_k)_{k\leq N}=(Z_k,S_k)_{k\leq N}$. This quantization consists in finding an optimally designed discretization of the process to provide for each step $k$ the best possible approximation of $\Theta_k$ by a random variable $\widehat{\Theta}_k$ which state space has a finite and fixed number of points. Here, \emph{optimal} means that the distance between $\Theta_k$ and $\widehat{\Theta}_k$ in a suitably chosen $L^{p}$ norm is minimal. For details on the quantization methods, we mainly refer to \cite{pages04} but the interested reader can also consult \cite{bally03}, \cite{bally05} and the references therein.\\

More precisely, consider $X$ an $\R^{q}$-valued random variable such that $\|X\|_{p}<\infty$ and let $K$ be a fixed integer. The optimal $L_{p}$-quantization of the random variable $X$ consists in finding the best possible $L_{p}$-approximation of $X$ by a random vector $\widehat{X}$ taking at most $K$ values: $\widehat{X}\in \{x^{1},\ldots,x^{K}\}$.
This procedure consists in the following two steps:
\begin{enumerate}
\item Find a finite weighted grid $\Gamma\subset \mathbb{R}^q$ with $\Gamma= \{x^{1},\ldots,x^{K}\}$.
\item Set $\widehat{X}=\widehat{X}^{\Gamma}$ where $\widehat{X}^{\Gamma}=proj_{\Gamma}(X)$ with $proj_{\Gamma}$ denotes the closest neighbour projection on $\Gamma$.
\end{enumerate}

\noindent
The asymptotic properties of the $L_{p}$-quantization are given by the following result, see e.g. \cite{pages04}.
\begin{theorem}
\label{theore}
If $\mathbb{E}[|X|^{p+\eta}]<+\infty$ for some $\eta>0$ then one has
\begin{eqnarray*}
\lim_{K\rightarrow \infty} K^{p/q} \min_{|\Gamma|\leq K} \| X-\widehat{X}^{\Gamma}\|^{p}_{p}& =& J_{p,q} 
\int |h|^{q/(q+p)}(u)du,
\end{eqnarray*}
where the law of $X$ is $P_{X}(du)=h(u) \lambda_{q}(du)+\nu$ with $\nu\perp\lambda_{d}$, $J_{p,q}$ a constant and $ \lambda_{q}$ the Lebesgue measure in $\mathbb{R}^{q}$.
\end{theorem}
Remark that $X$ needs to have finite moments up to the order $p+\eta$ to ensure the above convergence.
In this work, we used the CLVQ quantization algorithm described in \cite{bally03}, Section 3.\\

There exists a similar procedure for the optimal quantization of a Markov chain $\{X_{k}\}_{k\in \N}$. There are two approaches to provide the quantized approximation of a Markov chain. The first one, based on the quantization at each time $k$ of the random variable $X_{k}$ is called the \textit{marginal quantization}. The second one that enhances the preservation of the Markov property is called \textit{Markovian quantization}. Remark that for the latter, the quantized Markov process is not homogeneous. These two methods are described in details in \cite[section 3]{pages04}. In this work, we used the marginal quantization approach for simplicity reasons.\\

The quantization algorithm provides for each time step $0\leq k\leq N$ a finite grid $\Gamma_k$ of $E\times\R^+$ as well as the transition matrices $(\widehat{Q}_k)_{0\leq k\leq N-1}$ from $\Gamma_k$ to $\Gamma_{k+1}$. Let~$p~\geq~1$ such that for all $k\leq N$, $Z_k$ and $S_k$ have finite moments at least up to order~$p$ and let $proj_{\Gamma_{k}}$ be the closest-neighbor projection from $E\times\R^+$ onto $\Gamma_k$ (for the distance associated to norm~$p$). The quantized process $(\widehat{\Theta}_k)_{k\leq N}=(\widehat{Z}_k,\widehat{S}_k)_{k\leq N}$ takes values for each $k$ in the finite grid $\Gamma_k$ of $E\times\R^+$ and is defined by
\begin{equation}(\widehat{Z}_k,\widehat{S}_k)=proj_{\Gamma_{k}}(Z_k,S_k).\label{Z_n hat Z_n-mesurable}\end{equation}

Moreover, we also denote respectively $\Gamma_k^{Z}$ and $\Gamma_k^{S}$ the projections of $\Gamma_k$ on $E$ and $\R^{+}$.\\

Some important remarks must be made concerning the quantization. On the one hand, the optimal quantization has nice convergence properties stated by Theorem \ref{theore}. Indeed, the $L^p$-quantization error $\|\Theta_k-\widehat{\Theta}_k\|_p$ goes to zero when the number of points in the grids goes to infinity. However, on the other hand, the Markov property is not maintained by the algorithm and the quantized process is generally not markovian. Although the quantized process can be easily transformed into a Markov chain (see \cite{pages04}), this chain will not be homogeneous. It must be pointed out that the quantized process $(\widehat{\Theta}_k)_{k\in\N}$ depends on the starting point $\Theta_0$ of the process.\\

In practice, we begin with the computation of the quantization grids which merely requires to be able to simulate the process. This step is quite time-consuming, especially when the number of points in the quantization grids is large. However, the grids are only computed once and for all and may be stored off-line. What is more, they only depend on the dynamics of the process, not on the cost functions $l$ and $c$. Hence, the same grids may be used to compute different expectations of functionals as long as they are related to the same process. Our schemes are then based on the following simple idea: we replace the process by its quantized approximation within the operator $G$. The approximation is thus obtained in a very simple way since the quantized process has finite state space.

\subsection{Approximation of the expectation and rate of convergence} \label{DefRecQuantif}\label{section-theo}

We now use the quantization of the process $(\Theta_{k})_{k\leq N}=(Z_{k},S_{k})_{k\leq N}$. In order to approximate the random variables $(V_{k})_{k\leq N}$, we introduce a quantized version of the operator $G$. Notice that the quantized process is no longer an homogeneous Markov chain so that we have different operators for each time step $k$. Their definitions naturally stem from a remark made in the previous section: recall that for all $k\in \{1,...,N\}$ and $x\in E$,
\begin{align*}
Gw(x)=\E\left[F(Z_{k-1},S_{k})+w(Z_{k})\big|Z_{k-1}=x\right]
\end{align*}
\begin{definition} For all $k\in \{1,...,N\}$, $w\in B(\Gamma_{k}^{Z})$ and $z\in \Gamma_{k-1}^{Z}$, let
$$\widehat{G}_{k}w(z)=\E\left[F(z,\widehat{S}_{k})+w(\widehat{Z}_{k})\big|\widehat{Z}_{k-1}=z\right],$$
we then introduce the functions $(\widehat{v}_{k})_{0\leq k\leq N}$:
$$\left\{\begin{array}{lll}
\widehat{v}_{N}(z)&=0,&\text{ for all } z\in \Gamma_N^{Z},\\
\widehat{v}_{k}(z)&=\widehat{G}_{k+1}\widehat{v}_{k+1}(z),&\text{ for all } k\in \{0,...,N-1\} \text{ and } z\in \Gamma_{k}^{Z}.
\end{array}\right.$$
Eventually, for all $k\in \{0,...,N\}$, let
$$\widehat{V}_{k}=\widehat{v}_{k}(\widehat{Z}_{k}).$$
\end{definition}
\begin{remarque}
The conditional expectation in $\widehat{G}_{k}w(z)$ is a finite sum. Thus, the numerical computation of the sequence $(\widehat{V}_{k})_{k}$ will be easily performed as soon as the quantized process $(\widehat{\Theta}_k)_{k\leq N}$ has been obtained.
\end{remarque}
\begin{remarque} We have assumed that $Z_{0}=x$ a.s. Thus, the quantization algorithm provides that $\widehat{Z}_{0}=x$ a.s. too. Consequently, the random variable $\widehat{V}_{0}=\widehat{v}_{0}(\widehat{Z}_{0})$ is, in fact, deterministic.

\end{remarque}

The following theorem states the convergence of $\widehat{V}_{0}$ towards $V_{0}=J_{N}^{A}(l,c)(x)$ and provides a bound for the rate of convergence.

\begin{theorem}\label{TheoConvExpect} For all $k\in \{0,...,N\}$, one has $v_{k}\in\Lc(E)$. Moreover, the approximation error satisfies: 
$$|J_{N}^{A}(l,c)(x)-\widehat{V}_{0}|\leq \varepsilon_{N}(l,c,X,A)$$
where
\begin{multline*}
\varepsilon_{N}(l,c,X,A)= \sum_{k=0}^{N-1}\Big(2[v_{k+1}]\|Z_{k+1}-\widehat{Z}_{k+1}\|_p \\+\big(2[v_{k}]+[F]_1\big)\|Z_{k}-\widehat{Z}_{k}\|_p+[F]_{2}\|S_{k+1}-\widehat{S}_{k+1}\|_p\Big)+\frac{NC_{c}C_{\lambda}}{A}
\end{multline*}
with
$$\begin{array}{lll}
&[F]_{1}&=C_{t^{*}}[l]_{1}+[c]_{*}+A[t]_{*}C_{c},\\
&[F]_{2}&=C_{l}+AC_{c}.\\
&C_{v_{n}}&\leq n \big(C_{t^*}C_{l} + C_{c}\big),\\
&[v_{n}]_1&\leq e^{C_{t^*}C_\lambda}\left(K(A,v_{n-1})+n C_{t^*}[\lambda]_1\Big(C_{t^*}C_{l} + C_{c}\Big)\right)+C_{\ts}[l]_{1},\\
&[v_{n}]_2&\leq e^{C_{t^*}C_\lambda}\Big(C_{\ts}C_{l}C_{\lambda}+2C_{l}+C_\lambda C_{c}+(2n-1)C_{\lambda} \big(C_{t^*}C_{l} + C_{c}\big)\Big)+C_{l},\\
&[v_{n}]_*&\leq [v_{n}]_1+[t^*][v_{n}]_2.\\
&[v_{n}]&\leq K(A,v_{n-1}),
\end{array}$$
and for all $w\in\Lc(E)$, $K(A,w)=E_1+E_{2}A+E_3[w]_1+E_4C_w+[Q][w]_*$ where eventually
$$\begin{array}{ll}
E_1=&2[l]_{1}C_{\ts}+C_{l}\big([\ts]+2C_{\ts}^{2}[\lambda]_{1}\big)+[c]_{*}\big(1+C_{\ts}C_{\lambda}\big)\\
&+C_{c}\big(2[\lambda]_{1}C_{\ts}+C_{\lambda}C_{\ts}^{2}[\lambda]_{1}+2[\ts]C_{\lambda}\big),\\
E_2=&C_{c}C_{\ts}C_{\lambda}[\ts],\\
E_3=&\big(1+C_{\ts}C_{\lambda}\big)[Q],\\
E_4=&2C_{\lambda}[\ts]+C_{\ts}[\lambda]_{1}\big(2+C_{\ts}C_{\lambda}\big).
\end{array}$$
\end{theorem}

\paragraph{The choice of $A$}\label{choixA}Proposition~\ref{prop-JAconv} suggests that $A$ should be as large as possible. However, the constants $[F]_{1}$, $[F]_{2}$ and $[v_{n}]$ that appear in the bound of the approximation error proposed by the above theorem~\ref{TheoConvExpect} grow linearly with $A$. Thus, in order to control this error, it is necessary that the order of magnitude of the quantization error $\|\Theta_{k}-\widehat{\Theta}_{k}\|_{p}$ be at most $\frac{1}{A}$.\\

The convergence of the approximation scheme can be derived from theorem~\ref{TheoConvExpect}. Indeed, on the one hand, one must remind that $V_{0}=J_{N}^{A}(l,c)(x)$ is the expectation we intended to approximate and on the other hand, $\|\Theta_k-\widehat{\Theta}_k\|_p$ may become arbitrarily small when the number of points in the quantization grids goes to infinity (see e.g.~\cite{pages04}). An outline of the proof is presented in Appendix~\ref{proof-TheoConvExpect}.

\section{Time depending functionals}\label{section-transformation}

We now turn to the main contribution of this paper and present two generalizations of the previous problem. On the one hand, we will consider time depending functionals of the form
$$\E_{x}\left[\int_{0}^{T_{N}}l(X_{t},t)dt+\sum_{j=1}^{N}c(X_{T^{-}_{j}},T_{j})\1_{\{X_{T^{-}_{j}}\in \partial E\}}\right]$$
where $l$ and $c$ are Lipschitz continuous functions. On the other hand, we wish to replace the random time horizon $T_{N}$ by a deterministic one denoted $t_{f}$ i.e.
$$\E_{x}\left[\int_{0}^{t_{f}}l(X_{t},t)dt+\sum_{T_{j}\leq t_{f}}c(X_{T^{-}_{j}},T_{j})\1_{\{X_{T^{-}_{j}}\in \partial E\}}\right].$$
We will reason as follows. As it is suggested by M.H.A. Davis in \cite{davis93}, we will introduce a transformation $(\widetilde{X}_{t})_{t\geq 0}$ of the initial process $(X_{t})_{t\geq 0}$ by including the time variable into the state space i.e. $(\widetilde{X}_{t})=(X_{t},t)$. Indeed, we will see that both the expectation of the time depending functional and the one with deterministic time horizon are no other than expectations of time invariant functionals for the time augmented process $(\widetilde{X}_{t})_{t\geq 0}$. We therefore intend to apply the previously exposed approximation scheme to this new PDMP. However, it is far from obvious that the Lipschitz continuity assumptions \ref{hyp-lip-lambda}, \ref{hyp-lip-Q} and \ref{hyp-lip-t*} still hold for this new process.

Thus, the rest of this section is organized as follows. First, we will recall the precise definition of the time augmented process and prove that it satisfies the Lipschitz continuity assumptions required by our approximation scheme. Then, we will see that the time depending functional case corresponds to a time invariant functional for the new transformed process and may therefore be obtained thanks to the earlier method. Eventually, we will consider the deterministic time horizon problem that features an additional hurdle namely the presence of non Lipschitz continuous indicator functions.

\subsection{The time augmented process}\label{section-time-augemented}

M.H.A. Davis suggests, in \cite{davis93}, section 31, that the case of the time dependent functionals may be treated by introducing the time variable within the state space. Thus, it will be possible to apply our previous numerical method to the time augmented process. However, and this is what we discuss in this section, it is necessary to check whether the Lipschitz continuity assumptions still hold. We first recall the definition of the time augmented process given by M.H.A. Davis.

\begin{definition}
Introduce the new state space
$$\widetilde{E}=E\times \R^{+}$$
equipped with the norm defined by: for all $\xi=(x,t)$, $\xi'=(x',t')\in \widetilde{E}$, let
\begin{equation}\label{def-dist-tilde}
|\xi-\xi'|=|x-x'|+|t-t'|
\end{equation}
where the norm on $E$ is given by \eqref{def-dist}. On this state space, we define the process
$$\widetilde{X}_{t}=(X_{t},t).$$
\end{definition}
The local characteristics of the PDMP $(\widetilde{X}_{t})_{t\geq 0}$, denoted $(\widetilde{\lambda},\widetilde{Q},\widetilde{\Phi})$ are given for all $\xi=(x,t)\in \widetilde{E}$ by

\begin{align*}
\left\{\begin{array}{ll}
\widetilde{\lambda}(\xi)=\lambda(x),&\\
\widetilde{\Phi}\big(\xi,s\big)=\big(\Phi(x,s),t+s\big) &\text{ for $s\leq t^{*}(x)$,}\\
\widetilde{Q}\big(\xi,A\times \{t\}\big)=Q\big(x,A\big) &\text{ for all $A\in \mathcal{B}(E)$.}
\end{array}\right.
\end{align*}
Moreover, we naturally define for all $\xi=(x,t)\in \widetilde{E}$
$$\widetilde{t}^{*}(\xi)=\inf\{s>0 \text{ such that } \widetilde{\Phi}(\xi,s)\in\partial \widetilde{E}\}= t^{*}(x)$$

Clearly, Assumptions \ref{hyp-lip-lambda} and \ref{hyp-lip-t*} still hold with $[\widetilde{\lambda}]_{1}=[\lambda]_{1}$ and $[\widetilde{\ts}]=[\ts]$. However, proving assumption~\ref{hyp-lip-Q} is more intricate. We start with the following lemma.

\begin{lemme}\label{lemme-wt} Let $u,t\geq 0$ and $w\in \Lc^{u}(\widetilde{E})$. Denote $w_{t}$ the function of $B(E)$ defined by $w_{t}=w(\cdot,t)$. One has then $w_{t}\in \Lc^{t\wedge u}(E)$ with $[w_{t}]_{1}^{E,t\wedge u}\leq [w]_{1}^{\widetilde{E},u}$, $[w_{t}]_{2}^{E,t\wedge u}\leq [w]_{1}^{\widetilde{E},u}+[w]_{2}^{\widetilde{E},u}$, and $[w_{t}]_{*}^{E,t\wedge u}\leq (1+[\ts])[w]_{*}^{\widetilde{E},u}$.
\end{lemme}
\demo
Let $u,t\geq 0$ and $w\in \Lc^{u}(\widetilde{E})$. On the one hand, for $x,x'\in E$ and $s\leq t^{*}(x)\wedge t^{*}(x')\wedge t\wedge u$, one has
$$\left|w_{t}(\Phi(x,s))-w_{t}(\Phi(x',s))\right|=\left|w\Big(\widetilde{\Phi}\big((x,t-s),s\big)\Big)-w\Big(\widetilde{\Phi}\big((x',t-s),s\big)\Big)\right|.$$
We now use the fact that $w\in \Lc^{u}(\widetilde{E})$ which yields since $s\leq u$
$$\left|w_{t}(\Phi(x,s))-w_{t}(\Phi(x',s))\right|\leq[w]^{\widetilde{E},u}_{1}\big|(x,t-s)-(x',t-s)|=[w]^{\widetilde{E},u}_{1}\big|x-x'\big|.$$
Hence, $[w_{t}]^{E,t\wedge u}_{1}\leq[w]^{\widetilde{E},u}_{1}$ and similarly one obtains $[w_{t}]^{E,t\wedge u}_{2}\leq [w]_{1}^{\widetilde{E},u}+[w]_{2}^{\widetilde{E},u}$.\\

On the other hand, for $x,x'\in E$ such that $t^{*}(x)\vee t^{*}(x')\leq t\wedge u$, one~has
\begin{eqnarray*}\lefteqn{\left|w_{t}(\Phi(x,t^{*}(x)))-w_{t}(\Phi(x',t^{*}(x')))\right|}\\
&=&\left|w\Big(\widetilde{\Phi}\big((x,t-t^{*}(x)),t^{*}(x)\big)\Big)-w\Big(\widetilde{\Phi}\big((x',t-t^{*}(x')),t^{*}(x')\big)\Big)\right|\\
&=&\Big|w\Big(\widetilde{\Phi}\big((x,t-t^{*}(x)),\widetilde{t}^{*}(x,t-t^{*}(x))\big)\Big)\\
&& -w\Big(\widetilde{\Phi}\big((x',t-t^{*}(x')),\widetilde{t}^{*}(x',t-t^{*}(x'))\big)\Big)\Big|
\end{eqnarray*}
moreover since $w\in \Lc^{u}(\widetilde{E})$ and $\widetilde{t}^{*}(x,t-t^{*}(x))\vee \widetilde{t}^{*}(x',t-t^{*}(x'))\leq u$ one has
\begin{eqnarray*}\left|w_{t}(\Phi(x,t^{*}(x)))-w_{t}(\Phi(x',t^{*}(x')))\right|\leq[w]^{\widetilde{E},u}_{*}\big|(x,t-t^{*}(x))-(x',t-t^{*}(x'))\big|
\end{eqnarray*}
We conclude thanks to the Lipschitz continuity assumption~\ref{hyp-lip-t*} on $\ts$ providing
$\big|(x,t-t^{*}(x))-(x',t-t^{*}(x'))\big| \leq \big(1+[t^{*}]\big)\big|x-x'\big|$. Eventually, one has $[w_{t}]^{E,t\wedge u}_{*}\leq[w]^{\widetilde{E},u}_{*}\big(1+[t^{*}]\big)$ and $w_{t}\in \Lc^{t\wedge u}(E)$.
\findemo\\

The following proposition proves that Assumption~\ref{hyp-lip-Q} holds for the time augmented process $(\widetilde X)_{t\geq 0}$.

\begin{proposition}For all $w\in \Lc^{u}(\widetilde{E})$, one has
\begin{enumerate}
\item{for all $\xi$, $\xi'\in \widetilde{E}$ and $s\in [0,\widetilde{t}^*(\xi)\wedge \widetilde{t}^*(\xi')\wedge u]$,
$$\left|\widetilde{Q}w\big(\widetilde{\Phi}\big(\xi,s\big)\big)-\widetilde{Q}w\big(\widetilde{\Phi}\big(\xi',s\big)\big)\right| \leq ([Q]\vee 1)[w]^{\widetilde{E},u}_{1}\big|\xi-\xi'\big|,$$}
\item{for all $\xi$, $\xi'\in \widetilde{E}$ such that $\widetilde{t}^*(\xi)\vee \widetilde{t}^*(\xi')\leq u$,
\begin{multline*}
\left|\widetilde{Q}w\big(\widetilde{\Phi}\big(\xi,\widetilde{t}^{*}(\xi)\big)\big)-\widetilde{Q}w\big(\widetilde{\Phi}\big(\xi',\widetilde{t}^{*}(\xi')\big)\big)\right|\\\leq\big([Q]\vee 1\big)\big(1+[t^{*}]\big)\big([w]^{\widetilde{E},u}_{*}+[w]^{\widetilde{E},u}_{1}\big)|\xi-\xi'|,
\end{multline*}}
\end{enumerate}
in other words, Assumption~\ref{hyp-lip-Q} is satisfied with $[\widetilde{Q}]= \big([Q]\vee 1\big)\big(1+[t^{*}]\big)$.
\end{proposition}
\demo As in the previous lemma, for all $t\geq 0$, we will denote $w_{t}$ the function of $B(E)$ defined by $w_{t}=w(\cdot,t)$.
For $\xi=(x,t)\in \widetilde{E}$ and $w\in\Lc^{u}(\widetilde{E})$, one has, by definition of $\widetilde{Q}$,
\begin{equation}\label{eq-Qtilde-Q}
\widetilde{Q}w(\xi)=\int_{\xi'\in\widetilde{E}}w(\xi')\widetilde{Q}\big((x,t),d\xi'\big)
=\int_{z\in E}w(z,t)Q\big(x,dz\big)
=Qw_{t}(x).
\end{equation}

We may now check the regularity assumption on $\widetilde{Q}$. Let $\xi=(x,t)$ and $\xi'=(x',t')\in \widetilde{E}$. Let $s\in [0;\widetilde{t}^{*}(\xi)\wedge\widetilde{t}^{*}(\xi')\wedge u]$. Thanks to the definition of $\widetilde{\Phi}$ and equation \eqref{eq-Qtilde-Q} one has
\begin{eqnarray*}
\lefteqn{\left|\widetilde{Q}w\big(\widetilde{\Phi}\big(\xi,s\big)\big)-\widetilde{Q}w\big(\widetilde{\Phi}\big(\xi',s\big)\big)\right|}\\
&=&\left|\widetilde{Q}w\big(\Phi(x,s),t+s\big)-\widetilde{Q}w\big(\Phi(x',s),t'+s\big)\right|\\
&=&\left|Qw_{t+s}\big(\Phi(x,s)\big)-Qw_{t'+s}\big(\Phi(x',s)\big)\right|
\end{eqnarray*}
We split it into the sum of two differences
\begin{multline*}
\left|Qw_{t+s}\big(\Phi(x,s)\big)-Qw_{t'+s}\big(\Phi(x',s)\big)\right|\\
\leq \left|Qw_{t+s}\big(\Phi(x,s)\big)-Qw_{t+s}\big(\Phi(x',s)\big)\right|+\left|Q(w_{t+s}-w_{t'+s})\big(\Phi(x',s)\big)\right|.
\end{multline*}
On the one hand, we recall that thanks to lemma \ref{lemme-wt}, $w_{t+s}\in\Lc^{(t+s)\wedge u}(E)$ so that, since $s\leq (t+s)\wedge u$, we may use the Lipschitz continuity assumption~\ref{hyp-lip-Q} on $Q$ and the first term is bounded as follows 
$$\left|Qw_{t+s}\big(\Phi(x,s)\big)-Qw_{t+s}\big(\Phi(x',s)\big)\right|\leq[Q][w_{t+s}]^{E,(t+s)\wedge u}_{1}|x-x'|.$$ 
Moreover, lemma \ref{lemme-wt} also provides $[w_{t+s}]_{1}^{E,(t+s)\wedge u} \leq [w]^{\widetilde{E},u}_{1}$. On the other hand, and more basically, the second term in the above equation satisfies 
$$\left|Q(w_{t+s}-w_{t'+s})\big(\Phi(x',s)\big)\right|\leq [w]^{\widetilde{E},u}_{1}|t-t'|.$$ 
Eventually, one has
\begin{equation*}
\left|\widetilde{Q}w\big(\widetilde{\Phi}\big(\xi,s\big)\big)-\widetilde{Q}w\big(\widetilde{\Phi}\big(\xi',s\big)\big)\right|
\leq ([Q]\vee 1)[w]^{\widetilde{E},u}_{1}\big|\xi-\xi'\big|.
\end{equation*}

We now reason similarly to bound $\left|\widetilde{Q}w\big(\widetilde{\Phi}\big(\xi,\widetilde{t}^{*}(\xi)\big)\big)-\widetilde{Q}w\big(\widetilde{\Phi}\big(\xi',\widetilde{t}^{*}(\xi')\big)\big)\right|$ where $\xi=(x,t)$ and $\xi'=(x',t')\in \widetilde{E}$ are such that $\widetilde{t}^*(\xi)\vee \widetilde{t}^*(\xi')\leq u$. Equation~\eqref{eq-Qtilde-Q} yields
\begin{multline*}
\left|\widetilde{Q}w\big(\widetilde{\Phi}\big(\xi,\widetilde{t}^{*}(\xi)\big)\big)-\widetilde{Q}w\big(\widetilde{\Phi}\big(\xi',\widetilde{t}^{*}(\xi')\big)\big)\right|\\
=\left|Qw_{t+t^{*}(x)}\big(\Phi(x,t^{*}(x))\big)-Qw_{t'+t^{*}(x')}\big(\Phi(x',t^{*}(x'))\big)\right|
\end{multline*}
that we now spilt into
\begin{eqnarray*}
\lefteqn{\left|Qw_{t+t^{*}(x)}\big(\Phi(x,t^{*}(x))\big)-Qw_{t'+t^{*}(x')}\big(\Phi(x',t^{*}(x'))\big)\right|}\\
&\leq& \left|Qw_{t+t^{*}(x)}\big(\Phi(x,t^{*}(x))\big)-Qw_{t+t^{*}(x)}\big(\Phi(x',t^{*}(x'))\big)\right|\\
&&+\left|Q(w_{t+t^{*}(x)}-Qw_{t'+t^{*}(x')})\big(\Phi(x',t^{*}(x'))\big)\right|.
\end{eqnarray*}
Thanks to lemma \ref{lemme-wt}, $w_{t+t^{*}(x)}\in \Lc^{(t+t^{*}(x))\wedge u}(E)$. Moreover, we assume, without loss of generality that $\ts(x)\geq \ts(x')$ so that $\ts(x)\vee \ts(x') \leq (t+t^{*}(x))\wedge u$. Therefore, the first term in the above equation is bounded, thanks to the Lipschitz continuity assumption~\ref{hyp-lip-Q} on $Q$ and lemma \ref{lemme-wt}, by $[Q]\big((1+[t^{*}])[w]^{\widetilde{E}, u}_{*}+[w]^{\widetilde{E}, u}_{1}\big)|x-x'|$. More basically, the second term is bounded by $[w]^{\widetilde{E}, u}_{1}|t-t'+t^{*}(x)-t^{*}(x')|\leq [w]^{\widetilde{E}, u}_{1}(|t-t'|+[t^{*}]|x-x'|)$. Eventually, one has
\begin{align*}
&\left|\widetilde{Q}w\big(\widetilde{\Phi}\big(\xi,\widetilde{t}^{*}(\xi)\big)\big)-\widetilde{Q}w\big(\widetilde{\Phi}\big(\xi',\widetilde{t}^{*}(\xi')\big)\big)\right|\\
\leq&[Q](1+[t^{*}])[w]^{\widetilde{E}, u}_{*}|x-x'|+[w]^{\widetilde{E}, u}_{1}\big([Q]|x-x'|+|t-t'|+[\ts]|x-x'|\big)\\
\leq&\big([Q]\vee 1\big)\big(1+[t^{*}]\big)\big([w]^{\widetilde{E}, u}_{*}+[w]^{\widetilde{E}, u}_{1}\big)|\xi-\xi'|.
\end{align*}
Hence the result.
\findemo\\

Consequently, we may apply our numerical method to the time augmented process $(\widetilde{X}_{t})_{t\geq 0}$. In other words, for $l\in \Lc(\widetilde{E})$, $c\in\Lc(\partial \widetilde{E})$ and $\xi\in \widetilde{E}$, our approximation scheme may be used to compute
\begin{equation}\label{Jtilde}
\widetilde{J}_{N}(l,c)(\xi)=\E_{\xi}\left[\int_{0}^{T_{N}}l(\widetilde{X}_{t})dt+\sum_{j=1}^{N}c(\widetilde{X}_{T^{-}_{j}})\1_{\{\widetilde{X}_{T^{-}_{j}}\in \partial \widetilde{E}\}}\right].
\end{equation}
We will now see that the time depending functional and the deterministic time horizon problems boil down to computing such quantities $\widetilde{J}_{N}(l,c)(\xi)$ for suitably chosen functions $l$ and $c$. 

\subsection{Lipschitz continuous cost functions}

We first consider the time depending functional problem with Lipschitz continuous cost functions. Thus, let then $l\in \Lc(\widetilde{E})$, $c\in\Lc(\partial \widetilde{E})$ and $x\in E$, we wish to compute
$$\E_{x}\left[\int_{0}^{T_{N}}l(X_{t},t)dt+\sum_{j=1}^{N}c(X_{T^{-}_{j}},T_{j})\1_{\{X_{T^{-}_{j}}\in \partial E\}}\right].$$

It is straightforward to notice that this quantity may be expressed using the time augmented process starting from the point~$\xi_{0}=(x,0)$. Indeed, one has
$$\widetilde{J}_{N}(l,c)(\xi_{0})=\E_{x}\left[\int_{0}^{T_{N}}l(X_{t},t)dt+\sum_{j=1}^{N}c(X_{T^{-}_{j}},T_{j})\1_{\{X_{T^{-}_{j}}\in \partial E\}}\right]$$
where $\widetilde{J}_{N}(l,c)(\xi_{0})$ is given by equation \eqref{Jtilde}. Although they are time depending, the cost functions $l$ and $c$ are seen, in the left-hand side term, as time invariant functions of the time augmented process. The expectation of the time depending functional is therefore obtained by computing the expectation of a time invariant functional for the transformed PDMP thanks to the approximation scheme described in Section~\ref{section-approximation}. This is what expresses the following theorem, which proof stems from the previous discussion.

\begin{theorem} Let $l\in\Lc(\widetilde{E})$ and $c\in\Lc(\partial \widetilde{E})$ and apply the approximation scheme described in Section~\ref{section-approximation} to the time augmented process $(\widetilde{X}_{t})_{t\geq0}$, one has then
$$
\left|\E_{x}\left[\int_{0}^{T_{N}}l(X_{t},t)dt+\sum_{j=1}^{N}c(X_{T^{-}_{j}},T_{j})\1_{\{X_{T^{-}_{j}}\in \partial E\}}\right] - \widehat{V}_{0}\right|
 \leq \varepsilon_{N}(l,c,\widetilde{X},A).
$$
where we denoted $\varepsilon_{N}(l,c,\widetilde{X},A)$ the bound of the approximation error provided by Theorem~\ref{TheoConvExpect} when our approximation scheme is applied with cost functions $l$ and $c$ to the time augmented process $(\widetilde{X}_{t})_{t\geq0}$.
\end{theorem}

\begin{remarque}\label{rq-epsilon-Xtilde}The quantity $\varepsilon_{N}(l,c,\widetilde{X},A)$ is computed with respect to the process $(\widetilde{X}_{t})_{t\geq0}$ instead of $(X_{t})_{t\geq0}$ as presented in Theorem~\ref{TheoConvExpect} so that one has
\begin{multline*}
\varepsilon_{N}(l,c,\widetilde{X},A)=\sum_{k=0}^{N-1}\Big(2[v_{k+1}]^{\widetilde{E}}\|\widetilde{Z}_{k+1}-\widehat{\widetilde{Z}}_{k+1}\|_p \\+\big(2[v_{k}]^{\widetilde{E}}+[F]'_1+[F]''_1A\big)\|\widetilde{Z}_{k}-\widehat{\widetilde{Z}}_{k}\|_p\\
+\big([F]'_{2}+A[F]''_{2}\big)\|\widetilde{S}_{k+1}-\widehat{\widetilde{S}}_{k+1}\|_p\Big)+\frac{NC_{c}C_{\lambda}}{A}.
\end{multline*}
where $(\widetilde{Z}_{k},\widetilde{S}_{k})_{k\in\N}$ denotes the sequence of the post-jump locations and the inter-jump times of the time augmented process $(\widetilde{X}_{t})_{t\geq0}$ and with
$$\begin{array}{lll}
&[F]'_{1}=C_{t^{*}}[l]^{\widetilde{E}}_{1}+[c]^{\widetilde{E}}_{*},\\
&[F]''_{1}=[\ts]C_{c},\\
&[F]'_{2}=C_{l},\\
&[F]''_{2}=C_{c},\\
&C_{v_{n}}\leq n \big(C_{t^*}C_{l} + C_{c}\big),\\
&[v_{n}]^{\widetilde{E}}_1\leq e^{C_{t^*}C_\lambda}\left(\widetilde{K}(A,v_{n-1})+n C_{t^*}[\lambda]_1\Big(C_{t^*}C_{l} + C_{c}\Big)\right)+C_{\ts}[l]^{\widetilde{E}}_{1},\\
&[v_{n}]^{\widetilde{E}}_2\leq e^{C_{t^*}C_\lambda}\Big(C_{\ts}C_{l}C_{\lambda}+2C_{l}+C_\lambda C_{c}+(2n-1)C_{\lambda} \big(C_{t^*}C_{l} + C_{c}\big)\Big)+C_{l},\\
&[v_{n}]^{\widetilde{E}}_*\leq [v_{n}]^{\widetilde{E}}_1+[t^*][v_{n}]^{\widetilde{E}}_2.\\
&[v_{n}]^{\widetilde{E}}\leq \widetilde{K}(A,v_{n-1}),
\end{array}$$
and for all $w\in\Lc(E)$, $\widetilde{K}(A,w)=\widetilde{E}_1+E_{2}A+\widetilde{E}_3[w]^{\widetilde{E}}_1+E_4C_w+[\widetilde{Q}][w]^{\widetilde{E}}_*$ where eventually
$$\begin{array}{ll}
[\widetilde{Q}]=& \big([Q]\vee 1\big)\big(1+[t^{*}]\big),\\
\widetilde{E}_1=&2[l]^{\widetilde{E}}_{1}C_{\ts}+C_{l}\big([\ts]+2C_{\ts}^{2}[\lambda]_{1}\big)+[c]^{\widetilde{E}}_{*}\big(1+C_{\ts}C_{\lambda}\big)\\
&+C_{c}\big(2[\lambda]_{1}C_{\ts}+C_{\lambda}C_{\ts}^{2}[\lambda]_{1}+2[\ts]C_{\lambda}\big),\\
E_2=&C_{c}C_{\ts}C_{\lambda}[\ts],\\
\widetilde{E}_3=&\big(1+C_{\ts}C_{\lambda}\big)[\widetilde Q],\\
E_4=&2C_{\lambda}[\ts]+C_{\ts}[\lambda]_{1}\big(2+C_{\ts}C_{\lambda}\big).
\end{array}$$
\end{remarque}

\subsection{Deterministic time horizon}\label{section-deterministic-horizon}

In the context of applications, it seems relevant to consider a deterministic time horizon $t_{f}$. For instance, one may want to estimate a mean cost over a given period no matter how many jumps occur during this period. Actually, we will choose a time horizon of the form $t_{f}\wedge T_{N}$ with $N$ large enough to ensure the $N$-th jump will occur after time $t_{f}$ with a high probability i.e. that $\PP_{x}\big(T_{N}<t_{f}\big)$ be close to zero. For a discussion concerning the choice of such $N$, and in particular a theoretical bound of the probability $\PP_{x}\big(T_{N}<t_{f}\big)$, we refer to \cite{brandejsky10}. Simply notice that in practice, this probability may be estimated through Monte-Carlo simulations. We thus intend to approximate the following quantity for $l\in \Lc(\widetilde{E})$, $c\in\Lc(\partial \widetilde{E})$ and $x\in E$:
\begin{align*}
&\E_{x}\Big[\int_{0}^{T_{N}\wedge t_{f}}l(X_{t},t)dt+\sum_{T_{j}\leq t_{f}}c(X_{T^{-}_{j}},T_{j})\1_{\{X_{T^{-}_{j}}\in \partial E\}}\Big]\\
=&\E_{x}\Big[\int_{0}^{T_{N}}l(X_{t},t)\1_{\{t\leq t_{f}\}}dt+\sum_{j=1}^{N}c(X_{T^{-}_{j}},T_{j})\1_{\{X_{T^{-}_{j}}\in \partial E\}}\1_{\{T_{j}\leq t_{f}\}}\Big]
\end{align*}
The natural approach would consist in killing the process at time $t_{f}$ as M.H.A. Davis suggests in \cite{davis93}, section 31, and applying our method to the new process. However, the killed process will not necessarily fulfill our Lipschitz continuity assumptions because of the discontinuity introduced at time $t_{f}$.\\
A second idea would then be to use the previous results, to consider the time augmented process, and to define $\widetilde{l}(x,t)=l(x,t)\1_{\{t\leq t_{f}\}}$ and $\widetilde{c}(x,t)=c(x,t)\1_{\{t\leq t_{f}\}}$. However, a similar problem appears. Indeed, such functions $\widetilde{l}$ and $\widetilde{c}$ are not Lipschitz continuous and our numerical method requires this assumption. In the rest of this section, we will see how to overcome this drawback. On the one hand, we prove that the Lipschitz continuity condition on $l$ may be relaxed so that our numerical method may be used directly to approximate $\widetilde{J}_{N}(\widetilde{l},c)$ for any $c\in\Lc(\partial \widetilde{E})$. On the other hand, in the general case, we will deal with the non Lipschitz continuity of $\widetilde{c}$ by bounding it between two Lipschitz continuous functions.

\subsubsection{Direct estimation of the running cost term}

Let us explain how the Lipschitz continuity condition on the running cost function may be relaxed so that Theorem~\ref{TheoConvExpect}, stating the convergence of our approximation scheme, remains true when the running cost function is $\widetilde{l}(x,t)=l(x,t)\1_{\{t\leq t_{f}\}}$ with $l\in\Lc(\widetilde{E})$ and the boundary jump cost function is $c\in\Lc(\partial \widetilde{E})$ (although with slightly different constants in the bound of the convergence rate). Indeed, the running cost function $\widetilde l$ appears inside an integral that will have a regularizing effect allowing us to derive the required Lipschitz property of the functional in spite of the discontinuity of $\widetilde l$. Details are provided in Appendix~\ref{section-appendix-relax}.\\

Consequently, our approximation scheme may be used to compute $\widetilde{J}_{N}(\widetilde{l},c)(\xi)$ for any $\xi\in \widetilde{E}$. We recall that $\widetilde{J}_{N}$ is defined by \eqref{Jtilde} and that for all $x\in E$, one has
$$\widetilde{J}_{N}(\widetilde{l},c)(x,0)=\E_{x}\left[\int_{0}^{T_{N}\wedge t_{f}}l(X_{t},t)dt+\sum_{j=1}^{N}c(X_{T^{-}_{j}},T_{j})\1_{\{X_{T^{-}_{j}}\in \partial E\}}\right].$$
We now turn to the indicator function $\1_{\{T_{j}\leq t_{f}\}}$ required within the boundary jump cost term.

\subsubsection{Bounds of the boundary jump cost term}

We explained how the Lipschitz continuity condition on $l$ may be relaxed. However, when it comes to $c$, this condition cannot be avoided and our numerical method cannot be used directly with $\widetilde{c}(x,t)=c(x,t)\1_{\{t\leq t_{f}\}}$. We overcome this drawback by using Lipschitz continuous approximations of the indicator function. Indeed, for $B>0$, we introduce the real-valued functions $\underline{u}_{B}$ and $\overline{u}_{B}$ defined on $\R$ by
\begin{align*}
\underline{u}_{B}(t)&=\left\{\begin{array}{ll}
1 &\text{if $t< t_{f}-\frac{1}{B}$,}\\
-B(t-t_{f}) &\text{if $t_{f}-\frac{1}{B}\leq t <t_{f}$,}\\
0 &\text{if $t_{f}\leq t$,}\\
\end{array}\right.\\
\overline{u}_{B}(t)&=\left\{\begin{array}{ll}
1 &\text{if $t< t_{f}$,}\\
-B(t-t_{f})+1 &\text{if $t_{f}\leq t <t_{f}+\frac{1}{B}$,}\\
0 &\text{if $t_{f}+\frac{1}{B}\leq t$.}
\end{array}\right.
\end{align*}

\noindent The following lemma is straightforward.
\begin{lemme} For all $t\geq 0$, $\lim_{B\rightarrow +\infty}\underline{u}_{B}(t)=\1_{[0;t_{f}[}(t)$ and $\lim_{B\rightarrow +\infty}\overline{u}_{B}(t)=\1_{[0;t_{f}]}(t)$. Furthermore, for all $B>0$, $\underline{u}_{B}$ and $\overline{u}_{B}$ are Lipschitz continuous with Lipschitz constant $B$. Eventually, one has $\big|\underline{u}_{B}-\1_{[0;t_{f}]}\big|\leq 1$, $\big|\overline{u}_{B}-\1_{[0;t_{f}]}\big|\leq 1$ and $$\underline{u}_{B}\leq \1_{[0;t_{f}]} \leq \overline{u}_{B}.$$

\end{lemme}
\noindent Thus, define for $l\in\Lc(\widetilde{E})$ \begin{equation}\label{def-ltilde}\widetilde{l}(x,t)=l(x,t)\1_{\{t\leq t_{f}\}}\end{equation}
and for $c\in\Lc(\partial \widetilde{E})$ and for all $B>0$,

\begin{equation}\label{def-cb}\underline{c}_{B}(x,t)=c(x,t)\underline{u}_{B}(t) \qquad \text{and} \qquad \overline{c}_{B}(x,t)=c(x,t)\overline{u}_{B}(t).\end{equation}
We now check that these functions satisfy our Lipschitz continuity conditions.
\begin{proposition}\label{prop-lip-cb} The functions $\underline{c}_{B}$ and $\overline{c}_{B}$ belong to $\Lc(\partial \widetilde{E})$ with $[\underline{c}_{B}]_{*},[\overline{c}_{B}]_{*}\leq [c]_{*}+BC_{c}(1\vee[t^{*}])$.
\end{proposition}
\demo We prove the result for $\underline{c}_{B}$, the other case being similar. For all $\xi=(x,t),\xi'=(x',t') \in \widetilde{E}$, one has
\begin{eqnarray*}
\lefteqn{\big|\underline{c}_{B}\big(\widetilde{\Phi}(\xi,\ts(\xi))\big)-\underline{c}_{B}\big(\widetilde{\Phi}(\xi',\ts(\xi'))\big)\big|}\\
&=&\big|c\big(\widetilde{\Phi}(\xi,\widetilde{\ts}(\xi))\big)\underline{u}_{B}(t+\widetilde{\ts}(\xi))-c\big(\widetilde{\Phi}(\xi',\widetilde{\ts}(\xi'))\big)\underline{u}_{B}(t'+\widetilde{\ts}(\xi'))\big|\\
&\leq& [c]_{*}|\xi-\xi'|+C_{c}\big|\underline{u}_{B}(t+\widetilde{\ts}(\xi))-\underline{u}_{B}(t'+\widetilde{\ts}(\xi'))\big|\\
&\leq &[c]_{*}|\xi-\xi'|+C_{c}B\big(|t-t'|+[\widetilde{\ts}]|x-x'|\big)\\
&\leq &\big([c]_{*}+C_{c}B(1\vee[\ts])\big)|\xi-\xi'|.
\end{eqnarray*}
Hence the result.
\findemo\\

Therefore, the functions $\underline{c}_{B}$ and $\overline{c}_{B}$ are acceptable boundary jump cost functions and we may bound the deterministic horizon expectation by
\begin{multline*}
\widetilde{J}_{N}(\widetilde{l},\underline{c}_{B})(x,0)\leq\E_{x}\left[\int_{0}^{T_{N}}l(X_{t})\1_{\{t\leq t_{f}\}}dt+\sum_{j=1}^{N}c(X_{T^{-}_{j}})\1_{\{X_{T^{-}_{j}}\in \partial E\}}\1_{\{T_{j}\leq t_{f}\}}\right]\\ \leq \widetilde{J}_{N}(\widetilde{l},\overline{c}_{B})(x,0).
\end{multline*}
The following proposition provides the convergence of the bounds.
\begin{proposition} \label{prop-conv-deterministic-horizon}For all $x\in E$, one has
\begin{multline*}
\lim_{B\rightarrow +\infty}\widetilde{J}_{N}(\widetilde{l},\underline{c}_{B})(x,0)=\lim_{B\rightarrow +\infty}\widetilde{J}_{N}(\widetilde{l},\overline{c}_{B})(x,0)\\
=\E_{x}\left[\int_{0}^{T_{N}\wedge t_{f}}l(X_{t},t)dt+\sum_{j=1}^{N}c(X_{T^{-}_{j}},T_{j})\1_{\{X_{T^{-}_{j}}\in \partial E\}}\1_{\{T_{j}\leq t_{f}\}}\right].
\end{multline*}
The above convergence holds for every $t_{f}>0$ in the case of $\widetilde{J}_{N}(\widetilde{l},\overline{c}_{B})(x,0)$ but only for almost every $t_{f}>0$ with respect to the Lebesgue measure on $\R$ in the case of $\widetilde{J}_{N}(\widetilde{l},\underline{c}_{B})(x,0)$.
\end{proposition}
\demo Let $x\in E$. We first consider $\widetilde{J}_{N}(\widetilde{l},\overline{c}_{B})(x,0)$. 
\begin{align*}
\Big|\E_{x}&\left[\sum_{j=1}^{N}c(X_{T^{-}_{j}},T_{j})\1_{\{X_{T^{-}_{j}}\in \partial E\}}\1_{\{T_{j}\leq t_{f}\}}-\sum_{j=1}^{N}\overline{c}_{B}(X_{T^{-}_{j}},T_{j})\1_{\{X_{T^{-}_{j}}\in \partial E\}}\right]\Big|\\
&\leq\E_{x}\left[\sum_{j=1}^{N}\big|c(X_{T^{-}_{j}},T_{j})\big|\left|\1_{\{T_{j}\leq t_{f}\}}-\overline{u}_{B}(T_{j})\right|\right]\\
&\leq C_{c}\E_{x}\left[\sum_{j=1}^{N}\1_{\{t_{f} < T_{j} \leq t_{f}+\frac{1}{B}\}}\right]\\
&\leq C_{c}\sum_{j=1}^{N}\Big(\varphi_{j}(t_{f}+\frac{1}{B})-\varphi_{j}(t_{f})\Big)
\end{align*}
where $\varphi_{j}$ is the distribution function of $T_{j}$. For all $j\leq N$, $\lim_{B\rightarrow +\infty}\Big(\varphi_{j}(t_{f}+\frac{1}{B})-\varphi_{j}(t_{f})\Big)=0$ since $\varphi_{j}$ is right-continuous which shows the required convergence.\\

We now turn to the case of $\widetilde{J}_{N}(\widetilde{l},\underline{c}_{B})(x,0)$. Similar computations yields
\begin{align*}
\Big|\E_{x}&\left[\sum_{j=1}^{N}c(X_{T^{-}_{j}},T_{j})\1_{\{X_{T^{-}_{j}}\in \partial E\}}\1_{\{T_{j}\leq t_{f}\}}-\sum_{j=1}^{N}\underline{c}_{B}(X_{T^{-}_{j}},T_{j})\1_{\{X_{T^{-}_{j}}\in \partial E\}}\right]\Big|\\
&\leq C_{c}\sum_{j=1}^{N}\Big(\varphi_{j}(t_{f})-\varphi_{j}(t_{f}-\frac{1}{B})\Big).
\end{align*}
One cannot conclude as in the previous case since $\varphi_{j}$ is not necessary left-continuous. We therefore assume that $t_{f}$ is not an atom of any of the laws of the random variables $T_{j}$. Then, for all $j\leq N$, $\lim_{B\rightarrow +\infty}\Big(\varphi_{j}(t_{f})-\varphi_{j}(t_{f}-\frac{1}{B})\Big)=0$ and the result follows. Indeed, the set of the atoms of $T_{j}$ is at most countable so that the convergence holds for almost every $t_{f}$ w.r.t. the Lebesgue measure on $\R$.
\findemo

\subsubsection{Bounds in the general case}

The previous results show that the deterministic horizon expectation may be bounded by applying our numerical method with $\widetilde{l}$ and successively $\underline{c}_{B}$ and $\overline{c}_{B}$. In other words, we have shown the following theorem:

\begin{theorem} Let $l\in\Lc(\widetilde{E})$ and $c\in\Lc(\partial \widetilde{E})$. Let $(\underline{V}_{k,B})_{0\leq k\leq N}$ (respectively $(\overline{V}_{k,B})_{0\leq k \leq N}$) be the sequence of random variables $(V_{k})_{0\leq k\leq N}$ described in Section~\ref{section-approximation} when applying our approximation scheme to the time augmented process $(\widetilde{X}_{t})_{t\geq0}$ with cost functions $\widetilde{l}$ and $\underline{c}_{B}$ (resp. $\overline{c}_{B}$) defined by \eqref{def-ltilde} and \eqref{def-cb}. The bounds of the approximation error provided by Theorem~\ref{TheoConvExpect} are respectively denoted $\varepsilon_{N}(l,\underline{c}_{B},\widetilde{X},A,B)$ and $\varepsilon_{N}(l,\overline{c}_{B},\widetilde{X},A,B)$. One has then
\begin{multline*}
\underline{V}_{0,B}- \varepsilon_{N}(l,\underline{c}_{B},\widetilde{X},A,B)\\
\leq \E_{x}\left[\int_{0}^{T_{N}\wedge t_{f}}l(X_{t},t)dt+\sum_{j=1}^{N}c(X_{T^{-}_{j}},T_{j})\1_{\{X_{T^{-}_{j}}\in \partial E\}}\1_{\{T_{j}\leq t_{f}\}}\right]\\
\leq \overline{V}_{0,B}+ \varepsilon_{N}(l,\overline{c}_{B},\widetilde{X},A,B).
\end{multline*}
\end{theorem}

\begin{remarque}In the previous theorem, the quantity $\varepsilon_{N}(l,\underline{c}_{B},\widetilde{X},A,B)$ (and similarly $\varepsilon_{N}(l,\overline{c}_{B},\widetilde{X},A,B)$) is computed with respect to the process $(\widetilde{X}_{t})_{t\geq0}$ instead of $(X_{t})_{t\geq0}$ as presented in Theorem~\ref{TheoConvExpect} so that one has
\begin{multline*}
\varepsilon_{N}(l,\underline{c}_{B},\widetilde{X},A,B)=\sum_{k=0}^{N-1}\Big(2[v_{k+1}]^{\widetilde{E}}\|\widetilde{Z}_{k+1}-\widehat{\widetilde{Z}}_{k+1}\|_p \\+\big(2[v_{k}]^{\widetilde{E}}+[F]'_1+[F]''_1A+[F]'''_1B\big)\|\widetilde{Z}_{k}-\widehat{\widetilde{Z}}_{k}\|_p\\+\big([F]'_{2}+[F]''_{2}A\big)\|\widetilde{S}_{k+1}-\widehat{\widetilde{S}}_{k+1}\|_p\Big)+\frac{NC_{c}C_{\lambda}}{A}.
\end{multline*}
where $(\widetilde{Z}_{k},\widetilde{S}_{k})_{k\in\N}$ denotes the sequence of the post-jump locations and the inter-jump times of the time augmented process $(\widetilde{X}_{t})_{t\geq0}$ and with
$$\begin{array}{lll}
&[F]'''_{1}=C_{c}(1\vee[t^{*}]),\\
&[v_{n}]^{\widetilde{E}}_1\leq e^{C_{t^*}C_\lambda}\left(\widetilde{K}(A,B,v_{n-1})+n C_{t^*}[\lambda]_1\Big(C_{t^*}C_{l} + C_{c}\Big)\right)+C_{\ts}[l]^{\widetilde{E}}_{1},\\
&[v_{n}]^{\widetilde{E}}\leq \widetilde{K}(A,B,v_{n-1}),
\end{array}$$
and for all $w\in\Lc(E)$, $\widetilde{K}(A,B,w)=E_{1}'+E_{2}''B+E_{2}A+\widetilde{E}_3[w]^{\widetilde{E}}_1+E_4C_w+[\widetilde{Q}][w]^{\widetilde{E}}_*$ where eventually
$$\begin{array}{ll}
E'_1=&2[l]^{\widetilde{E}}_{1}C_{\ts}+C_{l}\big([\ts]+2C_{\ts}^{2}[\lambda]_{1}\big)+[c]^{\widetilde{E}}_{*}\big(1+C_{\ts}C_{\lambda}\big)\\
&+C_{c}\big(2[\lambda]_{1}C_{\ts}+C_{\lambda}C_{\ts}^{2}[\lambda]_{1}+2[\ts]C_{\lambda}\big),\\
E''_1=&C_{c}(1\vee[t^{*}])\big(1+C_{\ts}C_{\lambda}\big)\\
\end{array}$$
The other constants remain unchanged and we refer to remark \ref{rq-epsilon-Xtilde} for their precise expressions.
\end{remarque}

Furthermore, it is important to stress the fact that applying twice our numerical method does not increase significantly the computing time. Indeed, the computation of the quantization grids is, by far, the most costly step. These grids, that only depend on the dynamics of the process, may be stored off-line and used for the approximation of both bounds.

\paragraph{The choice of $B$.} We now discuss the choice of the parameter $B$, the discussion is quite similar to the one concerning the choice of $A$ in Section~\ref{section-theo}. proposition~\ref{prop-conv-deterministic-horizon} suggests that $B$ should be chosen as large as possible. However, choosing a large value for $B$ will lead to large Lipschitz constants that will decrease the sharpness of the bounds $\varepsilon_{N}(l,\underline{c}_{B},\widetilde{X})$ and $\varepsilon_{N}(l,\overline{c}_{B},\widetilde{X})$ for the approximation error provided by Theorem~\ref{TheoConvExpect}. Indeed, it is easy to check that $[v_{n}]$ grows linearly with $B$ (see the precise expressions of the Lipschitz constants above). Thus, in order to control the error proposed by Theorem~\ref{TheoConvExpect}, it is necessary that the order of magnitude of the quantization error $\|\Theta_{n}-\widehat{\Theta}_{n}\|_{p}$ be at most $\frac{1}{B}$.\\

\section{Numerical results}\label{section-result}

\subsection{A repair workshop model}

\indent We now present a repair workshop model adapted from \cite{davis93}, section 21.\\

In a factory, a machine produces goods which daily value is $r(x)$ where $x\in[0;1]$ represents a parameter of evolution of the machine, a setting chosen by the operator. For instance, $x$ may be some load or some pace imposed on the machine. This machine, initially working, may break down with age-dependent hazard rate $\lambda(t)$ and is then sent to the workshop for reparation. Besides, the direction of the factory has decided that, whenever the machine has worked for a whole year without requiring reparation, it is sent to the workshop for maintenance. The daily cost of such a maintenance is $q(x)$ while the daily cost of a reparation is $p(x)$, with reasonably $p(x)>q(x)$. We assume that after a reparation or a maintenance, that both last a fixed time $s$, the machine is totally repaired and is not worn down.

\begin{figure}[h!]
\begin{center}
\includegraphics[width=9cm]{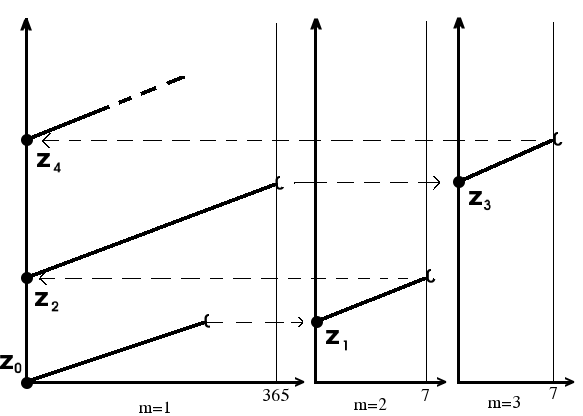}
\caption{An example of trajectory. The process starts from the point $Z_{0}$ in mode $m=1$ (machine in service). The machine may be sent to the workshop for repairs (m=2) or maintenance (m=3).}
\label{model-repar}
\end{center}
\end{figure}

We therefore consider three modes: the machine is working ($m=1$), being repaired ($m=2$), undergoing maintenance ($m=3$). The state of the process at time $t$  will be denoted $X_{t}=(m_{t},\zeta_{t},t)$ where $\zeta_{t}$ is the time since the last change of mode (this component is required since the hazard rate $\lambda$ is age-dependent). The state space is $E=\big(\{1\}\times [0;365]\times \R^{+} \big)\union\big(\{2\}\times [0;s]\times \R^{+} \big)\union\big(\{3\}\times [0;s]\times \R^{+} \big)$. In each mode, the flow is $\Phi_{m}\big((\zeta,t),u\big)=(\zeta+u,t+u)$. Concerning the transition kernel, from the previous discussion, one may notice for instance that from the point $(1,\zeta,t)$, the process can jump to the point $(2,0,t)$ if $\zeta<365$ and the jump is forced to $(3,0,t)$ if $\zeta=365$. Figure \ref{model-repar} presents the state space and an example of trajectory of the process.\\

Our aim is to find the value of the setting $x$ that maximizes the expected total benefits $B(x)$ i.e. the discounted value (the interest rate is denoted $\rho$) of production minus the costs related to maintenance and reparation over a period $t_{f}=5$ years:

$$B^{*}=\sup_{x\in[0;1]} B(x)$$
where
$$B(x)=\E_{(1,0,0)}\left[\int_{0}^{t_{f}}e^{-\rho t}\big(r(x)\1_{\{m_{t}=1\}}-p(x)\1_{\{m_{t}=2\}}-q(x)\1_{\{m_{t}=3\}}\big)dt\right].$$

We will use the following values $r(x)=x$, $p(x)=100x^{2}$, $q(x)=5$, $s=7$~days, $\rho=\frac{0.03}{365}$ and $\lambda$ represents a Weibull distribution with parameters $\alpha=2$ et $\beta=600$.

Our assumptions clearly hold so that we may run our numerical method. We first need to find $N\in \N$ such that $\PP_{(1,0,0)}(T_{N}<t_{f})$ be small. Monte-Carlo simulations lead to the value $N=18$. For a fixed $x\in[0;1]$, we will therefore compute $\widetilde{J}_{N}(\widetilde{l},0)(1,0,0)$ where $\widetilde{l}(m,\zeta,t)=e^{-\rho t}\big(r(x)\1_{\{m=1\}}-p(x)\1_{\{m=2\}}-q(x)\1_{\{m=3\}}\big)\1_{\{t\leq t_{f}\}}$. Eventually, notice that we could have chosen $r$, $p$ and $q$ slightly more generally by allowing them to be time-dependent.\\

It is important to stress the fact that, once the Markov chain associated to the process is quantized, we will be able to compute the approximation of $B(x)$ almost instantly for any $x\in[0;1]$ because the same grids are used for every computation. Thanks to this flexibility, we are able to draw the function $x\rightarrow B(x)$ and, thus, to solve the above optimization problem very easily. This is a very important advantage of our method. Indeed, if we computed $B(x)$ through standard methods such as Monte Carlo simulations, we would have to repeat the whole algorithm again and again for each value of $x$ and solving the optimization problem would be intractable.

The following figure represents the approximation of the function $B$ computed on a constant step grid of $[0;1]$ with step $10^{-2}$. This leads to the solution of the earlier optimization problem. Indeed, we obtain $B^{*}=B(x^{*})=537.84$ where $x^{*}=0.78$ is the value of the setting $x$ that maximizes the benefits of the factory.

\begin{figure}[h!]
\begin{center}
\includegraphics[width=9cm]{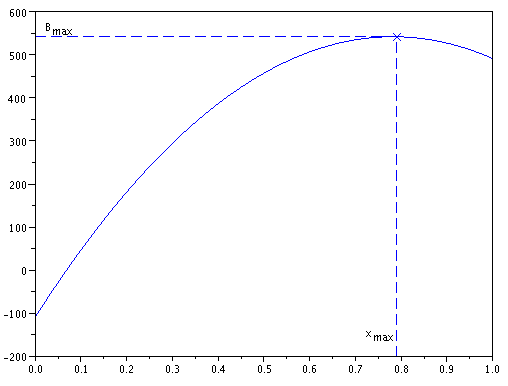}
\caption{The function $B$ drawn with 500 points in the quantization grids.}
\end{center}
\end{figure}

Let now $x=0.78$, the following table presents the values of $\widehat{V}_{N}$, which are the approximations of $B(x)$, for different number of points in the quantization grids. A reference value is obtained thanks to Monte-Carlo method ($10^{8}$ simulations) $B_{Monte-Carlo}=537.69$.

\begin{center}
\begin{tabular}{ccc}
\hline
Points in the quantization grids &$\widehat{V}_{N}$ & relative error to $537.69$\\
\hline
$20$ points& 542.14 & 0.83 \% \\
\hline
50 points & 539.57 &  0.35 \%\\
\hline
$100$ points& 538.24 &  0.10 \%\\
\hline
$500$ points& 537.84 & 0.03 \%\\
\hline
\end{tabular}
\end{center}

From a computational time point of view, we already explained that the computation of large quantization grids is, by far, the most costly step since it may take up to several hours whereas the approximation of the expectation that follows is then almost instantaneous. However, we may notice, in the above table, that grids containing only 50 points yield a quite accurate result with merely 0.35 \% error. Such grids only require a few minutes to be designed.

\begin{remarque} We already noticed that the same grids may serve several purposes. For instance, we may also have been interested in the computation of the mean time spent by the machine in the workshop by taking $l(m,\zeta,t)=\1_{\big\{m\in\{2;3\}\big\}}$.
\end{remarque}

\subsection{A corrosion model}\label{exemple-corrosion}

We consider here a corrosion model of an aluminum metallic structure. This example was provided by Astrium. It concerns a small structure within a strategic ballistic missile. The missile is stored successively in three different environments which are more or less corrosive. It is made to have potentially large storage durations. The requirement for security is very strong. The mechanical stress exerted on the structure depends in part on its thickness. A loss of thickness will cause an over-constraint and therefore increase a risk of rupture. It is thus crucial to study the evolution of the thickness of the structure over time.

Let us describe more precisely the usage profile of the missile. It is stored successively in three different environments, the workshop (denoted $m=1$), the submarine in operation ($m=2$) and the submarine in dry-dock ($m=3$). This is because the structure must be equipped and used in a given order. Then it goes back to the workshop and so on. The missile stays in each environment during a random duration with exponential distribution. Its parameter $\lambda_{m}$ depends on the environment. The degradation law for the thickness loss then depends on the environment through two parameters, a deterministic transition period $\eta_{m}$ and a random corrosion rate $\rho$ uniformly distributed within a given range. Typically, the workshop and dry-dock are the most corrosive environments but the time spent in operation is more important. The randomness of the corrosion rate accounts for small variations and uncertainties in the corrosiveness of each environment.

In each environment $m\in\{1;2;3\}$, the evolution over time of the thickness loss $d_{m}$ satisfies:
\begin{equation}\label{corrosion-def-d}
d_m(\rho,s)=\rho\left(s+\eta_m\left(e^{-\frac{s}{2\eta_m}}-1\right)\right).
\end{equation}
Table \ref{parametres cor} gives the values of the different parameters.\\

\begin{table}[h]
\begin{center}
\begin{tabular}{ccccc}
\multicolumn{1}{c}{}&&\text{environment 1}&\text{environment 2}&\text{environment 3}\\
\hline
$\lambda_m$ &($\text{h}^{-1}$)& $(17520)^{-1}$ & $(131400)^{-1}$ & $(8760)^{-1}$\\
\hline
$\eta_m$ &(h)& 30000 & 200000& 40000 \\
\hline
$\rho$ &(mm.$\text{h}^{-1}$)& $[10^{-6},10^{-5}]$ &  $[10^{-7},10^{-6}]$ &  $[10^{-6},10^{-5}]$ \\
\hline
\end{tabular}
\end{center}
\caption{Numerical values of the parameters of the corrosion model.}
\label{parametres cor}
\end{table}

Initially, the structure is in environment $m=1$ and the thickness loss is null. One draws the corrosion rate $\rho_{0}$ uniformly distributed in the interval $[10^{-6},10^{-5}]$ and the time of the first change of environment $T_{1}$ exponentially distributed with parameter $\lambda_{1}=(17520)^{-1}$ hours$^{-1}$. The corrosion starts according to Equation~\eqref{corrosion-def-d} so that, for all $0\leq t\leq T_{1}$, the loss of thickness is $d_{1}(\rho_{0},t)$. The structure then moves to environment 2 and the process restarts similarly: a new corrosion rate $\rho_{T_{1}}$ is drawn according to an uniform law on $[10^{-7},10^{-6}]$, the time of the second jump $T_{2}$ is drawn so that $T_{2}-T_{1}$ is exponentially distributed with parameter $\lambda_{2}=(131400)^{-1}$ hours$^{-1}$ and for $T_{1}\leq t\leq T_{2}$, the loss of thickness is $d_{1}(\rho_{0},T_{1})+d_{2}(\rho_{T_{1}},t-T_{1})$ and so on.

At each change of environment, a new corrosion rate $\rho$ is drawn according to a uniform law on the corresponding interval. The thickness loss, however, evolves continuously.

We are interested in computing the mean loss of thickness in environment 2 until a given time $t_{f}=18$ years.\\

\subsubsection*{Modelization by PDMP}

\paragraph{The state space $E$.}The loss of thickness will be modelized by a PDMP whose modes are the different environments. Let then $M=\{1,2,3\}$. The PDMP $(X_{t})_{t\geq 0}$ will contain the following components: the mode $m\in M$, the loss of thickness $d$, the time since the last jump $s$ (this is to ensure that the Markov property is satisfied), the corrosion rate $\rho$ and the time $t$ (since we consider the time-augmented process). Notice that clearly, one has always $s\leq t$ so that we reasonably consider the following state space:
$$E=\left\{(m,d,s,\rho,t)\in M\times \R^{+}\times \R^{+} \times [10^{-7};10^{-5}] \times \R^{+} \text{ such that } s\leq t\right\}.$$


\paragraph{The flow $\Phi$.}The flow is given for all $u\geq 0$ by
\begin{align*}
\Phi(\left(\begin{array}{c}m\\ d\\ s\\ \rho \\ t\end{array}\right),u)&=\left(\begin{array}{c}m\\ d + d_m(\rho,s+u) - d_m(\rho,s)\\ s+u \\ \rho \\ t+u\end{array}\right).
\end{align*}


\paragraph{The transition kernel $Q$.}Let us now study the jumps of this process. When the process jumps from a point $x=(m,d,s,\rho,t)\in E$, $m$ becomes $m+1$ modulo 3 (denoted $m+1 [3])$, $d$ and $t$ remain unchanged, $s$ becomes 0. Only $\rho$ is randomly drawn, according to a uniform law on an interval $[\rho_{min};\rho_{max}]$ that depends on the new mode. One has then for $w\in B(E)$, $x=(m,d,s,\rho,t)\in E$, and $u\geq 0$,
\begin{multline}\label{corrosion-Q}
Qw\big(\Phi(\left(\begin{array}{c}m\\ d\\ s\\ \rho \\ t\end{array}\right),u)\big)
=Qw\left(\begin{array}{c}m\\ d + d_m(\rho,s+u) - d_m(\rho,s)\\ s+u \\ \rho\\ t+u\end{array}\right)\\
=\frac{1}{\rho_{max}-\rho_{min}}\int_{\rho_{min}}^{\rho_{max}}w\left(\begin{array}{c}m+1[3]\\ d + d_m(\rho,s+u) - d_m(\rho,s)\\ 0 \\ \tilde\rho \\ t+u\end{array}\right)d\tilde\rho.
\end{multline}


\paragraph{The cost function $l$.}The function $l\in B(E)$ will be the cost function to compute the mean loss of thickness in mode 2. It is defined as follows: for all $x=(m,d,s,\rho,t)\in E$ and $u\geq 0$
\begin{equation}\label{corrosion-cost-l}
l(\Phi(x,u))=\rho\Big(1-\frac{1}{2}e^{-\frac{s+u}{2\eta_{m}}}\Big)\1_{\{m=2\}}=\frac{d}{du}\Big(d_{m}(\rho,s+u)\Big)\1_{\{m=2\}}.
\end{equation}
One then defines $\tilde{l}(\Phi(x,u))=l(\Phi(x,u))\1_{\{t+u\leq t_{f}\}}$ so that, one has
\begin{eqnarray*}
L(x,u)&=&\int_{0}^{u}\tilde{l}\big(\Phi(x,u')\big)du'\\
&=&\int_{0}^{u\wedge (t_{f}-t)^{+}}l\big(\Phi(x,u')\big)du'\\
&=&\Big(d_{m}(\rho,s+u\wedge (t_{f}-t)^{+})-d_{m}(\rho,s)\Big)\1_{\{m=2\}}
\end{eqnarray*}
that is indeed the thickness lost in mode $2$ from the point $x=(m,d,s,\rho,t)$ during a time $u\wedge (t_{f}-t)^{+}$.

\paragraph{The assumptions.}Assumption \ref{Tk_goes_to_infty} and \ref{hyp-lip-lambda} are clearly satisfied. Moreover, it is straightforward, from \eqref{corrosion-cost-l}, to check that $l\in \Lc(E)$ so that Assumption \ref{hyp-lip-lc} holds.\\

We now turn to Assumption \ref{hyp-lip-Q} and we will see that, although it does not hold for any function $w\in \Lc^{v}(E)$, it holds for a sufficiently big sub class of functions. We first need to make a remark. Recall that for all $x=(m,d,s,\rho,t)\in E$ and for all $k\in\{0,...,N\}$, one has $v_{N-k}(x)=\E_{x}\left[\int_{0}^{T_{k}}l\big(\Phi(x,u)\big)\1_{\{t+u\leq t_{f}\}}du\right]$. Therefore, for all $k\in\{0,...,N\}$ the function $v_{k}$ as well as the function $\tilde{l}$ satisfy the following condition:
\begin{equation}\label{w-null-apres-tf}
\text{for all }x=(m,d,s,\rho,t)\in E \text{ such that } t\geq t_{f}, \text{ one has } w(x)=0.
\end{equation}

The next step consists in proving that Assumption \ref{hyp-lip-Q}, although it is not satisfied for any function $w\in \Lc^{v}(E)$, holds for any function $w\in\Lc^{v}(E)$ that also satisfies condition \eqref{w-null-apres-tf}. This is done in Lemma \ref{lem-hyp-Q-cor} and it is sufficient because in the proof of the theorem that ensures the convergence of our approximation scheme, Assumption \ref{hyp-lip-Q} is only used with the functions $(v_{k})_{k\in\{0,...,N\}}$ that do satisfy condition \eqref{w-null-apres-tf}.
\begin{lemme}\label{lem-hyp-Q-cor} There exists $[Q]\in\R^{+}$ such that for all $v\geq 0$ and $w\in\Lc^{v}(E)$ that satisfies condition \eqref{w-null-apres-tf}, one has for all $x$, $x'\in E$ and $0\leq u\leq v$,
$$\left|Qw\big(\Phi(x,u)\big)-Qw\big(\Phi(x',u)\big)\right|\leq [Q][w]_{1}^{E,v}|x-x'|.$$
\end{lemme}
\demo Let $x=(m,d,s,\rho,t)$ and $x'=(m',d',s',\rho',t')\in E$ with for instance $t\leq t'$. First we may choose $m=m'$, otherwise, $|x-x'|=+\infty$ and there is nothing to prove. Now, we are facing three different cases:
\begin{itemize}
\item{if $t_{f}\leq t+u\leq t'+u$, then one has $Qw\big(\Phi(x,u)\big)=Qw\big(\Phi(x',u)\big)=0$ because $w$ satisfies condition \eqref{w-null-apres-tf} and there is nothing to prove.}
\item{if $t+u\leq t_{f}\leq t'+u$, notice that $Qw\big(\Phi(x',u)\big)=Qw\big(\Phi((m',d',s',\rho',t_{f}),u)\big)=0$ (this stems from condition \eqref{w-null-apres-tf}) so that we are reduced to the following case,}
\item{We assume from now on that $t+u\leq t'+u\leq t_{f}$. We now intend to bound $\left|Qw\big(\Phi(x,u)\big)-Qw\big(\Phi(x',u)\big)\right|$. It is clear from equation \eqref{corrosion-Q} that we only need to prove that the function $(\rho,s)\rightarrow d_{m}(\rho,s)$, defined by \eqref{corrosion-def-d}, is Lipschitz continuous w.r.t. both its variables on the set $[10^{-7};10^{-5}]\times [0;t_{f}]$. Indeed, we have $s\leq t$ and $s'\leq t'$ so that $s,s',s+u,s'+u\leq t_{f}$. Standard computations yield:
\begin{eqnarray*}
\left|d_{m}(\rho,s)-d_{m}(\rho',s')\right|&\leq& s|\rho-\rho'|+\frac{3}{2}\rho'|s-s'|\\
&\leq& t_{f}|\rho-\rho'|+\frac{3}{2}10^{-5}|s-s'|.
\end{eqnarray*}}
\end{itemize}
Hence the result.\findemo\\

Eventually, Assumption \ref{hyp-lip-t*} is not satisfied because in our corrosion model, one has $\ts(x)=+\infty$ for all $x\in E$. Besides, we may notice that the previous proof would have been more straightforward if $\ts$ had been bounded. Indeed in that case, we would have had $s,s',s+u,s'+u\leq C_{\ts}$ and the introduction of condition \eqref{w-null-apres-tf} would have been unnecessary. Nevertheless, we have been able to overcome the drawback of having $\ts$ non-bounded by noticing that somehow the deterministic time horizon $t_{f}$ plays the part of the missing $C_{\ts}$. This is the meaning of condition \eqref{w-null-apres-tf}: roughly speaking, we do not consider what happens beyond $t_{f}$.\\ 
More generally, we will now see that in our deterministic time horizon problem, the boundedness of $\ts$ may be dropped and our results remain true replacing $C_{\ts}$ by $t_{f}$. This is clear in the case of Proposition \ref{LipF} because the function $\tilde{l}$ satisfies the condition \eqref{w-null-apres-tf}. Proposition \ref{LipV_n} remains also true replacing $C_{\ts}$ by $t_{f}$. Indeed, on the one hand, it is clear that $L(x,u)\leq t_{f}C_{l}$. On the other hand, when computing $|v_{n}(\Phi(x,u))-v_{n}(\Phi(x',u'))|$, we are facing three different cases (as in the proof of Lemma \ref{lem-hyp-Q-cor}):
\begin{itemize}
\item{if $t_{f}\leq u\leq u'$, one has $v_{n}(\Phi(x,u))=v_{n}(\Phi(x',u'))=0$ (this stems from condition \eqref{w-null-apres-tf}),}
\item{if $u\leq t_{f}\leq u'$, one has $|v_{n}(\Phi(x,u))-v_{n}(\Phi(x',u'))|=|v_{n}(\Phi(x,u))-v_{n}(\Phi(x',t_{f}))|$ since $v_{n}(\Phi(x',u'))=v_{n}(\Phi(x',t_{f}))=0$ (condition \eqref{w-null-apres-tf} once again) so that we are reduced to the following case,}
\item{If $u\leq u'\leq t_{f}$, the computations remain unchanged and $t_{f}$ replaces $C_{\ts}$ as a bound for $u$ and $u'$.}
\end{itemize}

\subsubsection*{Numerical results}

\begin{figure}[h!]
\begin{center}
\includegraphics[width=10cm]{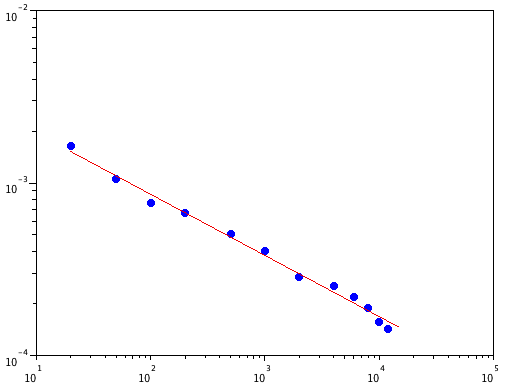}
\caption{Logarithm of the error when approximating the loss of thickness in environment 2 w.r.t. the logarithm of the number of points in the quantization grids. The empirical convergence rate, estimated through a regression model, is~-0.35.}
\label{err_cor_amb2}
\end{center}
\end{figure}

Table \ref{err_cor_table} presents the values of the loss of thickness in environment 2 obtained through our approximation scheme as well as a Monte Carlo approximation (obtained with $10^{8}$ simulations) and the relative errors of our values w.r.t. the Monte Carlo value. Figure \ref{err_cor_amb2} presents respectively the empirical convergence rate. The convergence rate, estimated through a regression model is -0.35. This is roughly the same order of magnitude as the rate of convergence of the optimal quantizer (see for instance \cite{pages04}) since here the dimension is 3 (indeed, $m$ is deterministic and $s=0$ immediately after a jump so that we only quantize the variables $\rho$, $d$ and $t$).

\begin{table}[h]
\begin{center}
\begin{tabular}{ccccccc}
\text{Quantization grids} &$\widehat V_{0}$&\text{error}\\
\hline
\text{20 points} &0.038386& 4.43 \%\\
\hline
\text{50 points}  &0.037804&2.85 \%\\
\hline
\text{100 points} &0.037525&2.09 \%\\
\hline
\text{200 points}  &0.037421&1.81 \%\\
\hline
\text{500 points}  &0.037264&1.38 \%\\
\hline
\text{1000 points} &0.037160&1.10 \%\\
\hline
\text{2000 points}  &0.037041&0.77 \%\\
\hline
\text{4000 points}  &0.037007&0.69 \%\\
\hline
\text{6000 points}  &0.036973&0.57 \%\\
\hline
\text{8000 points}  &0.036944&0.49 \%\\
\hline
\text{10000 points} &0.036911&0.40 \%\\
\hline
\text{12000 points} &0.036897&0.36 \%\\
\hline
\hline
\text{Monte Carlo} &0.036755 &\\
\hline
\end{tabular}
\end{center}
\caption{Approximation of the mean loss of thickness (in mm) in environment~2 for different numbers of points in the quantization grids and a Monte Carlo approximation ($10^{8}$ simulations).}
\label{err_cor_table}
\end{table}

Besides, Table \ref{cpu_time_table} presents the CPU time to compute the expectations from the quantization grids (computations are run with Matlab R2010b on a MacBook Pro 2.66 GHz i7 processor). It can be seen that, once the quantization grids are obtained, our approximation scheme is performed very quickly.

\begin{table}[h]
\begin{center}
\begin{tabular}{cc}
\text{Quantization grids}&\text{CPU time (second)}\\
\hline
\text{20 points} & 0.0059\\ 
\hline
\text{50 points} & 0.0085\\ 
\hline
\text{100 points} & 0.014\\ 
\hline
\text{200 points} & 0.034\\
\hline
\text{500 points} & 0.12\\ 
\hline
\text{1000 points} & 0.37\\
\hline
\text{2000 points} & 1.5\\
\hline
\text{4000 points} & 5.6\\
\hline
\text{6000 points} & 13\\
\hline
\text{8000 points} & 24 \\
\hline
\text{10000 points} & 35\\
\hline
\text{12000 points} & 54\\
\hline
\hline
\text{Monte Carlo ($10^{8}$ simulations)} & $\approx$16000\\
\hline
\end{tabular}
\end{center}
\caption{CPU time.}
\label{cpu_time_table}
\end{table}

\section{Conclusion}\label{section-conclusion}

We have presented an efficient and easy to implement numerical method to approximate expectations of functionals of piecewise-deterministic Markov processes. We proved the convergence of our algorithm with bounds for the rate of convergence.

Although our method concerns time invariant functionals, we proved that we are able to tackle time depending problems such as Lipschitz continuous time depending functionals or deterministic time horizon expectations. Indeed, we proved that, thanks to the introduction of the time augmented process, time depending problems may be seen, paradoxically, as special cases of the time invariant situation.

Our method is easy to implement because it merely requires to be able to simulate the process. Furthermore, although the computation of the quantization grids may be quite time-consuming, it may be performed preliminarily because the grids only depend on the dynamics of the process and not on the cost functions $l$ and $c$. Therefore, they may be stored off-line and serve several purposes. As it is illustrated by the examples presented in Section~\ref{section-result}, storing the grids provides to our approximation scheme efficiency and flexibility. Indeed, the computation of the expectation can be performed very quickly once the grids are available. Thus, if one decides for instance to modify the functional, the same grids may be used so that the new result is obtained very quickly. This flexibility is an important advantage over standard Monte-Carlo simulations.

\section*{Acknowledgements}

This work was supported by ARPEGE program of the French National Agency of Research (ANR),
project ''FAUTOCOES'', number ANR-09-SEGI-004. Besides, the authors gratefully acknowledge Astrium for its financial support.

\appendix

\section{Lipschitz continuity of $F$, $G$ and $v_{n}$}

The first lemma and the first proposition of this section present mainly the Lipschitz continuity of the functions $\delta^{A}$ and $F$. They are stated without proof because they are quite straightforward.
\begin{lemme}\label{Lip-delta}
The function $\delta^{A}$ is Lipschitz continuous w.r.t. both its variables i.e. for all $x$, $y\in E$ and $u$, $t\in \R$, one has
\begin{align*}
|\delta^{A}(x,t)-\delta^{A}(y,t)|&\leq A[\ts]|x-y|,\\
|\delta^{A}(x,t)-\delta^{A}(x,u)|&\leq A|t-u|,
\end{align*}
Moreover, one has for all $x\in E$ and $t,s\geq0$ such that $t+s\leq \tstar(x)$,
$$\delta^{A}(\Phi(x,s),t)=\delta^{A}(x,t+s).$$
\end{lemme}

\begin{proposition}\label{LipF} The function $F$, introduced in Definition \ref{def-L-C-F-G}, is Lipschitz continuous w.r.t. both its variables. For all $x$, $y\in E$ and $u,v\in [0;\ts(x)\wedge\ts(y)]$, one has

$$|F(x,u)-F(y,v)|\leq [F]_{1}|x-y|+[F]_{2}|u-v|$$
with
$$\begin{array}{lll}
&[F]_{1}&=C_{t^{*}}[l]_{1}+[c]_{*}+A[\ts]C_{c},\\
&[F]_{2}&=C_{l}+AC_{c}.
\end{array}$$
\end{proposition}

The two following lemmas are adapted from \cite{saporta10}, the second one being a special case of lemma A.1 from \cite{saporta10}. Thus, they are stated without proof.

\begin{lemme}\label{LipTech}For $h\in \mathbf{L_{c}}(E)$, $(x,y)\in E^2$, and $t\leq t^*(x)\wedge t^*(y)$
\begin{multline*}
\left|\int_t^{t^*(x)} h(\Phi(x,s))e^{-\Lambda(x,s)}ds-\int_t^{t^*(y)} h(\Phi(y,s))e^{-\Lambda(y,s)}ds\right|\\
\leq\Big(C_{t^*}[h]_1+\big(C_{t^*}^2[\lambda]_1+[t^*]\big)C_h\Big)|x-y|.
\end{multline*}
\end{lemme}

\begin{lemme}\label{LipTech-t*}For $h\in \mathbf{L_{c}}(\partial E)\union\mathbf{L_{c}}(E)$ and $x,y\in E$, one has
\begin{multline*}
\Big|e^{-\Lambda(x,t^*(x))}h\big(\Phi(x,t^*(x))\big)-e^{-\Lambda(y,t^*(y))}h\big(\Phi(y,t^*(y))\big)\Big|\\
\leq \Big([h]_*+C_h\big(C_{t^*}[\lambda]_1+[t^*]C_\lambda\big)\Big)|x-y|.
\end{multline*}
\end{lemme}

We now introduce a definition that will be convenient in the sequel. For $w \in \mathbf{L_{c}}(E)$, $x\in E$ and $t\in[0;t^{*}(x)]$, we define
\begin{align*}
G_{t}w(x)&=E_x\left[\left(F(x,S_1) + w(Z_{1})\right)\mathbbm{1}_{\{S_1\geq t\}}\right]\\
&=E_x\left[\left(L(x,S_1) + C(x,S_1) + w(Z_{1})\right)\mathbbm{1}_{\{S_1\geq t\}}\right].
\end{align*}
In particular, $G_{0}=G$. Since we know the law of $(Z_{1},S_{1})$, it can be shown that
\begin{equation}\label{Gdevelop}G_{t}w (x)=\Upsilon_{1}(x)+\Upsilon_{2}(x)+\Upsilon_{3}(x)+\Upsilon_{4}(x)+\Upsilon_{5}(x) \end{equation}
with
$$\begin{array}{ll}
\Upsilon_{1}(x)&=e^{-\Lambda(x,t)}\int_{0}^{t} l\circ\Phi(x,s)ds,\\
\Upsilon_{2}(x)&=\int_{t}^{t^*(x)} l\circ\Phi(x,s)e^{-\Lambda(x,s)}ds.\\
\Upsilon_{3}(x)&=c\circ\Phi(x,\ts(x))\int_{t}^{\ts(x)}\delta^{A}(x,s)\lambda\circ\Phi(x,s)e^{-\Lambda(x,s)}ds,\\
\Upsilon_{4}(x)&=\int_t^{t^*(x)}\big(\lambda Qw\big)\circ\Phi(x,s)e^{-\Lambda(x,s)}ds,\\
\Upsilon_{5}(x)&=e^{-\Lambda(x,t^*(x))}\big(Qw+c\big)\circ\Phi(x,t^*(x)).
\end{array}$$

\begin{proposition}\label{LipG}For $w \in \mathbf{L_{c}}(E)$, $(x,y)\in E^2$ and $t\in[0;t^{*}(x)\wedge t^{*}(y)]$, one has
$$\big|G_{t}w(x)-G_{t}w(y)\big|\leq K(A,w)|x-y|,$$ where $K(A,w)=E_1+E_{2}A+E_3[w]_1+E_4C_w+[Q][w]_*$ with
$$\begin{array}{ll}
E_1=&2[l]_{1}C_{\ts}+C_{l}\big([\ts]+2C_{\ts}^{2}[\lambda]_{1}\big)+[c]_{*}\big(1+C_{\ts}C_{\lambda}\big)\\
&+C_{c}\big(2[\lambda]_{1}C_{\ts}+C_{\lambda}C_{\ts}^{2}[\lambda]_{1}+2[\ts]C_{\lambda}\big),\\
E_2=&C_{c}C_{\ts}C_{\lambda}[\ts],\\
E_3=&\big(1+C_{\ts}C_{\lambda}\big)[Q],\\
E_4=&2C_{\lambda}[\ts]+C_{\ts}[\lambda]_{1}\big(2+C_{\ts}C_{\lambda}\big).
\end{array}$$
\end{proposition}

\demo Let $w \in \mathbf{L_{c}}(E)$, $(x,y)\in E^2$ and $t\in[0;t^{*}(x)\wedge t^{*}(y)]$. From equation~\eqref{Gdevelop}, we naturally split $|G_{t}w(x)-G_{y}w(y)|$ into the sum of five differences. \\

The first one is $|\Upsilon_{1}(x)-\Upsilon_{1}(y)|$ and is bounded by
\begin{align*}
|\Upsilon_{1}(x)-\Upsilon_{1}(y)|&\leq C_{\ts}C_{l}\Big|e^{-\Lambda(x,t)}-e^{-\Lambda(y,t)}\Big|+\int_{0}^{t} \Big(l\circ\Phi(x,s)-l\circ\Phi(y,s)\Big)ds\\
&\leq \Big(C_{\ts}^{2}C_{l}[\lambda]_{1} + C_{\ts}[l]_{1}\Big)|x-y|.
\end{align*}

The differences $|\Upsilon_{2}(x)-\Upsilon_{2}(y)|$ and $|\Upsilon_{4}(x)-\Upsilon_{4}(y)|$ may be bounded thanks to lemma \ref{LipTech} with successively $h=l$ and $h=\lambda Qw$. Notice that $C_{\lambda Qw}\leq C_\lambda C_w$ and $[\lambda Qw]_1\leq C_\lambda[Q][w]_1+C_w[\lambda]_1$. \\

For the difference of the $\Upsilon_{5}$ terms, we use lemma \ref{LipTech-t*} with $h=Qw+c$. Notice that $C_{Qw+c}\leq C_{w}+C_{c}$ and that $[Qw+c]_*\leq [Q]\big([w]_*+[w]_{1}\big)+[c]_{*}$.

Eventually, in order to bound $|\Upsilon_{3}(x)-\Upsilon_{3}(y)|$, we assume without loss of generality that $t^*(x)\leq t^*(y)$ and we have
\begin{eqnarray*}
\lefteqn{|\Upsilon_{3}(x)-\Upsilon_{3}(y)|}\\
&\leq& C_{c}\int_{t}^{\ts(x)}\Big|\delta^{A}(x,s)\lambda\circ\Phi(x,s)e^{-\Lambda(x,s)}-\delta^{A}(y,s)\lambda\circ\Phi(y,s)e^{-\Lambda(y,s)}\Big|ds\\ 
&&+C_{c}\int_{\ts(x)}^{\ts(y)}\Big|\delta^{A}(y,s)\lambda\circ\Phi(y,s)e^{-\Lambda(y,s)}\Big|ds+[c]_{*}C_{\ts}C_{\lambda}|x-y| \\
&\leq& C_{c}\int_{t}^{\ts(x)}\left(C_{\lambda}\big|\delta^{A}(x,s)-\delta^{A}(y,s)\big|+[\lambda]_{1}|x-y|+C_{\lambda}\big|e^{-\Lambda(x,s)}-e^{-\Lambda(y,s)}\big|\right)ds\\
&&+C_{c}[\ts]C_{\lambda}|x-y|+[c]_{*}C_{\ts}C_{\lambda}|x-y| \\
&\leq& \Big(C_{c}C_{\ts}\big(C_{\lambda}A[\ts]+[\lambda]_{1}+C_{\lambda}C_{\ts}[\lambda]_{1}\big)+C_{c}[\ts]C_{\lambda}+[c]_{*}C_{\ts}C_{\lambda}\Big)|x-y|.
\end{eqnarray*}
The result follows.
\findemo\\

The following lemma is stated without proof. It is indeed very close to Lemma 51.7 from \cite{davis93}.
\begin{lemme}\label{lemme-value-sur-le-flux}For all $x\in E$ and $t\in[0;t^{*}(x)]$, one has $$v_{n}(\Phi(x,t))=e^{\Lambda(x,t)}G_{t}v_{n-1}(x)-\int_{0}^{t} l\circ\Phi(x,s)ds.$$
\end{lemme}

\begin{proposition}\label{LipV_n}For all $n\in\{0,1,...,N\}$, $v_{n}\in \mathbf{L_{c}}(E)$ and one has
\begin{align*}
C_{v_{n}}&\leq n \big(C_{t^*}C_{l} + C_{c}\big),\\
[v_{n}]_1&\leq e^{C_{t^*}C_\lambda}\left(K(A,v_{n-1})+n C_{t^*}[\lambda]_1\Big(C_{t^*}C_{l} + C_{c}\Big)\right)+C_{\ts}[l]_{1},\\
[v_{n}]_2&\leq e^{C_{t^*}C_\lambda}\Big(C_{\ts}C_{l}C_{\lambda}+2C_{l}+C_\lambda C_{c}+(2n-1)C_{\lambda} \big(C_{t^*}C_{l} + C_{c}\big)\Big)+C_{l},\\
[v_{n}]_*&\leq [v_{n}]_1+[t^*][v_{n}]_2.\\
[v_{n}]&\leq K(A,v_{n-1}),
\end{align*}
\end{proposition}
\demo
Recall that for $x\in E$, one has from definition \ref{def-L-C-F-G}
$$v_{n}(x)=Gv_{n-1}(x)=\E_{x}\left[L(x,S_{1})\right]+\E_{x}\left[C(x,S_{1})\right]+\E_x[v_{n-1}(Z_1)].$$
Thus, $C_{v_{n}}\leq C_{t^*}C_{l} + C_{c}+C_{v_{n-1}}\leq n \big(C_{t^*}C_{l} + C_{c}\big)$ by induction.\\

Let us now turn to $[v_{n}]_1$. Lemma \ref{lemme-value-sur-le-flux} yields
\begin{eqnarray*}
\lefteqn{|v_{n}(\Phi(x,t))-v_{n}(\Phi(y,t))|}\\
&\leq&|e^{\Lambda(x,t)}G_{t}v_{n-1}(x)-e^{\Lambda(y,t)}G_{t}v_{n-1}(y)|+\int_{0}^{t} \big|l\circ\Phi(x,s)-l\circ\Phi(y,s)\big|ds\\
&\leq&e^{\Lambda(x,t)}|G_{t}v_{n-1}(x)-G_{t}v_{n-1}(y)|+|G_{t}v_{n-1}(y)||e^{\Lambda(x,t)}-e^{\Lambda(y,t)}|+C_{\ts}[l]_{1}|x-y|.
\end{eqnarray*}
The result follows using proposition~\ref{LipG} and noticing that
\begin{align*}
\Lambda(x,t)&\leq C_{t^*}C_\lambda,\\
|G_{t}v_{n-1}(y)|&\leq C_{t^*}C_{l} + C_{c}+C_{v_{n-1}}\leq n \big(C_{t^*}C_{l} + C_{c}\big),\\
|e^{\Lambda(x,t)}-e^{\Lambda(y,t)}|&\leq e^{C_{t^*}C_\lambda}C_{t^*}[\lambda]_1|x-y|.
\end{align*}
We now turn to $[v_{n}]_2$. For $x\in E$ and $s$, $t \in [0,t^*(x)]$ with $s\leq t$, one has
\begin{multline*}
|v_{n}(\Phi(x,t))-v_{n}(\Phi(x,s))|\\
\leq e^{\Lambda(x,t)}|G_{t}v_{n-1}(x)-G_{s}v_{n-1}(x)|+|G_{s}v_{n-1}(x)||e^{\Lambda(x,t)}-e^{\Lambda(x,s)}|+C_{l}|t-s|.
\end{multline*}
Moreover, from \eqref{Gdevelop}, one has
\begin{eqnarray*}
\lefteqn{|G_{t}v_{n-1}(x)-G_{s}v_{n-1}(x)|}\\
&\leq& E_x\left[\big|F(x,S_1) + v_{n-1} (Z_{1})\big|\mathbbm{1}_{\{s\leq S_1<t\}}\right]\\
&\leq& \Big|e^{-\Lambda(x,t)}\int_{0}^{t} l(\Phi(x,u))du-e^{-\Lambda(x,s)}\int_{0}^{s} l(\Phi(x,u))du\Big|\\
&&+\int_s^{t}\Big|l(\Phi(x,u))e^{-\Lambda(x,u)}\Big|du \\
&&+ \big|c\circ\Phi(x,\ts(x))\big|\int_{s}^{t}\Big|\delta^{A}(x,u)\lambda\circ\Phi(x,u)e^{-\Lambda(x,u)}\Big|du\\
&&+ \int_s^{t}\Big|\big(\lambda Qv_{n-1}\big)\circ\Phi(x,u)e^{-\Lambda(x,u)}\Big|du,\\
&\leq& \Big(C_{\ts}C_{l}|e^{-\Lambda(x,t)}-e^{-\Lambda(x,s)}|+C_{l}|t-s|\Big)+\Big(C_{l}|t-s|\Big)\\
&&+\Big(C_{c}C_{\lambda}|t-s|\Big)+\Big(C_{\lambda}C_{v_{n-1}}|t-s|\Big).
\end{eqnarray*}
and
$$|e^{\Lambda(x,t)}-e^{\Lambda(x,s)}|\leq e^{C_{t^*}C_\lambda}C_\lambda|t-s|.$$
Eventally, the bound for $[v_{n}]$ is a direct consequence from proposition~\ref{LipG}.
\findemo

\section{Relaxed assumption on the running cost function}\label{section-appendix-relax}

In this section, we consider the approximation applied to the time augmented process so that the local characteristics are $\widetilde{\Phi}$, $\widetilde{\lambda}$ and $\widetilde{Q}$ defined in Section~\ref{section-time-augemented}. Moreover, we consider a function $l\in \Lc(\widetilde{E})$ and we define $\widetilde{l}\in B(\widetilde{E})$ by
$$\text{for all $\xi=(x,t)\in \widetilde{E}$, }\widetilde{l}(\xi)=l(x,t)\1_{\{t\leq t_{f}\}}.$$
We intend to prove that the convergence of our approximation scheme, stated by Theorem~\ref{TheoConvExpect}, remains true if we choose $\widetilde{l}$ as the running cost function even though it does not fulfill the required Lipschitz conditions i.e. $\widetilde{l}\not\in \Lc(\widetilde{E})$. Indeed, the Lipschitz continuity of $l$ is used four times in the proof of the theorem, once in proposition~\ref{LipF}, twice in proposition~\ref{LipG} (when bounding the difference of the $\Upsilon_{1}$ terms and the one of the $\Upsilon_{2}$ ones) and once in proposition~\ref{LipV_n} (when bounding $[v_{n}]_{1}$). In each case, the Lipschitz continuity of the running cost function $l$ is used to bound a term of the form
\begin{equation}\label{form-ltilde1}
\int_{s}^{s'}\Big|\widetilde{l}\circ\widetilde{\Phi}\big(\xi,u\big)-\widetilde{l}\circ\widetilde{\Phi}\big(\xi',u\big)\Big|du
\end{equation}
for $\xi$, $\xi' \in \widetilde{E}$ and $s$, $s'\in [0;\widetilde{t}^{*}(\xi)\wedge \widetilde{t}^{*}(\xi')]$, or of the form
\begin{equation}\label{form-ltilde2}
\int_s^{\widetilde{t}^*(\xi)\wedge \widetilde{t}^{*}(\xi')} \left|\widetilde{l}\circ\widetilde{\Phi}\big(\xi,u\big)e^{-\widetilde{\Lambda}\big(\xi,u\big)}-\widetilde{l}\circ\widetilde{\Phi}\big(\xi',u\big)e^{-\widetilde{\Lambda}\big(\xi',u\big)}\right|du
\end{equation}
for $\xi$, $\xi' \in \widetilde{E}$ and $s\in [0;\widetilde{t}^{*}(\xi)\wedge \widetilde{t}^{*}(\xi')]$ and where we naturally denoted $\widetilde{\Lambda}(\xi,u)=\int_0^u \widetilde{\lambda}(\widetilde{\Phi}(\xi,v))dv$. Concerning this second form, equation \eqref{form-ltilde2}, notice that
\begin{eqnarray*}
\lefteqn{\int_s^{\widetilde{t}^*(\xi)\wedge \widetilde{t}^{*}(\xi')} \left|\widetilde{l}\circ\widetilde{\Phi}\big(\xi,u\big)e^{-\widetilde{\Lambda}\big(\xi,u\big)}-\widetilde{l}\circ\widetilde{\Phi}\big(\xi',u\big)e^{-\widetilde{\Lambda}\big(\xi',u\big)}\right|du}\\
&\leq& \int_s^{\widetilde{t}^*(\xi)\wedge \widetilde{t}^{*}(\xi')} \left|\widetilde{l}\circ\widetilde{\Phi}\big(\xi,u\big)-\widetilde{l}\circ\widetilde{\Phi}\big(\xi',u\big)\right|du\\
&&+ C_{l}\int_s^{\widetilde{t}^*(\xi)\wedge \widetilde{t}^{*}(\xi')} \left|e^{-\widetilde{\Lambda}\big(\xi,u\big)}-e^{-\widetilde{\Lambda}\big(\xi',u\big)}\right|du\\
&\leq& \int_s^{\widetilde{t}^*(\xi)\wedge \widetilde{t}^{*}(\xi')} \left|\widetilde{l}\circ\widetilde{\Phi}\big(\xi,u\big)-\widetilde{l}\circ\widetilde{\Phi}\big(\xi',u\big)\right|du + C_{l}C_{\ts}^{2}[\lambda]_{1}
\end{eqnarray*}
so that, to ensure that Theorem~\ref{TheoConvExpect} remains true with $\widetilde l$ as the running cost function, it is sufficient to be able to bound terms of the form \eqref{form-ltilde1}. This is done in the following lemma.

\begin{lemme}\label{lemme-Lip-ltilde}For $\xi=(x,t),\xi'=(x',t') \in \widetilde{E}$ and $s\in [0;\widetilde{t}^{*}(\xi)\wedge \widetilde{t}^{*}(\xi')]$, one has
$$\int_0^s\Big|\widetilde{l}\circ\widetilde{\Phi}\big(\xi,u\big)-\widetilde{l}\circ\widetilde{\Phi}\big(\xi',u\big)\Big|du\leq (C_{t^{*}}[l]_{1}+ C_{l})|\xi-\xi'|.$$
\end{lemme}
\demo
Let $\xi=(x,t),\xi'=(x',t') \in \widetilde{E}$ and $s\in [0;\widetilde{t}^{*}(\xi)\wedge \widetilde{t}^{*}(\xi')]$, one has
\begin{eqnarray*}
\lefteqn{\int_0^s\Big|\widetilde{l}\circ\widetilde{\Phi}\big(\xi,u\big)-\widetilde{l}\circ\widetilde{\Phi}\big(\xi',u\big)\Big|du}\\
&\leq&\int_0^s\Big|l\circ\widetilde{\Phi}\big(\xi,u\big)\1_{\{t+u\leq t_{f}\}}-l\circ\widetilde{\Phi}\big(\xi',u\big)\1_{\{t'+u\leq t_{f}\}}\Big|du\\
&\leq&\int_0^s\Big|l\circ\widetilde{\Phi}(\xi,u)-l\circ\widetilde{\Phi}(\xi',u)\Big|du+C_{l}\int_0^s\Big|\1_{\{t+u\leq t_{f}\}}-\1_{\{t'+u\leq t_{f}\}}\Big|du
\end{eqnarray*}
The left-hand side term is bounded by $C_{t^{*}}[l]_{1}|\xi-\xi'|$ since $l\in\Lc(\widetilde{E})$. For the right-hand side term, assume without loss of generality that $t\leq t'$, one has
$$\Big|\1_{\{t+u\leq t_{f}\}}-\1_{\{t'+u\leq t_{f}\}}\Big|=\Big|\1_{\{t-t_{f}\leq u\}}-\1_{\{t'-t_{f}\leq u\}}\Big|=\1_{\{t-t_{f}\leq u< t'-t_{f}\}}$$
so that the right-hand side term is bounded by $C_{l}|t-t'|\leq C_{l}|\xi-\xi'|$. The result follows.
\findemo\\

Eventually, Theorem~\ref{TheoConvExpect}, remains true if we choose $\widetilde{l}$ as the running cost function. One only needs to slightly modify the Lipschitz constants given in propositions \ref{LipF}, \ref{LipG} and \ref{LipV_n}. The terms $C_{t^{*}}[l]_{1}$ have to be replaced by $C_{t^{*}}[l]_{1}+ C_{l}$.

\section{Proof of Theorem~\ref{TheoConvExpect}}\label{proof-TheoConvExpect}

The Lipschitz continuity of the functions $v_{k}$ is proved by proposition~\ref{LipV_n}. Let now $A>0$ and first notice that 
$$|J_{N}(l,c)(x)-\widehat{V}_{0}|\leq |J_{N}(l,c)(x)-V_{0}|+|V_{0}-\widehat{V}_{0}|.$$
Proposition \ref{prop-JAconv} states that $|J_{N}(l,c)(x)-V_{0}|\leq \frac{NC_{c}C_{\lambda}}{A}$ since $V_{0}=J^{A}_{N}(l,c)(x)$. We now have to bound $|V_{0}-\widehat{V}_{0}|$.\\

Some of the arguments of the proof are similar to the ones used in Theorem~5.1 from \cite{saporta10}, thus we will not develop the details of the proof. Recall that $\|V_{N}-\widehat{V}_{N}\|_p=0$ and let $k\in\{0,...,N-1\}$. In order to bound the approximation error, let us split it into three terms $\|V_{k}-\widehat{V}_{k}\|_p\leq \Xi_1+\Xi_2+\Xi_3$ where
\begin{align*}\left\{\begin{array}{ll}
\Xi_1&=\|v_{k}(Z_{k})-v_{k}(\widehat{Z}_{k})\|_p,\\
\Xi_2&=\|G v_{k+1}(\widehat{Z}_{k})-\widehat{G}_{k+1} v_{k+1}(\widehat{Z}_{k})\|_p,\\
\Xi_3&=\|\widehat{G}_{k+1} v_{k+1}(\widehat{Z}_{k})-\widehat{G}_{k+1} \widehat{v}_{k+1}(\widehat{Z}_{k})\|_p.
\end{array}\right.\end{align*}
The theorem is then a direct consequence from the three following lemmas, stated without proof, that provide bounds for each of these three terms.
\begin{lemme}
The first term $\Xi_1$ is bounded by
$$\|v_{k}(Z_{k})-v_{k}(\widehat{Z}_{k})\|_p\leq [v_{k}]\|Z_{k}-\widehat{Z}_{k}\|_p.$$
\end{lemme}

\begin{lemme}The second term $\Xi_2$ is bounded by
\begin{multline*}\|G v_{k+1}(\widehat{Z}_{k})-\widehat{G}_{k+1} v_{k+1}(\widehat{Z}_{k})\|_p\\
\leq [v_{k+1}]\|Z_{k+1}-\widehat{Z}_{k+1}\|_p+\big([v_{k}]+[F]_1\big)\|Z_{k}-\widehat{Z}_{k}\|_p+[F]_{2}\|S_{k+1}-\widehat{S}_{k+1}\|_p.
\end{multline*}
\end{lemme}

\begin{lemme}
The third term $\Xi_3$ is bounded by
\begin{multline*}
\|\widehat{G}_{k+1} v_{k+1}(\widehat{Z}_{k})-\widehat{G}_{k+1} \widehat{v}_{k+1}(\widehat{Z}_{k})\|_p 
\leq [v_{k+1}]\|Z_{k+1}-\widehat{Z}_{k+1}\|_p + \|V_{k+1}-\widehat{V}_{k+1}\|_p.
\end{multline*}
\end{lemme}

\bibliographystyle{plain}
\bibliography{biblio}

\end{document}